\documentclass[12pt]{article}
\usepackage{amsmath,amsfonts,amssymb,amsthm}

\newcommand{\I}{{\bf 1}}
\newtheorem{proposition}{Proposition}[section]
\newtheorem{theorem}[proposition]{Theorem}

\newtheorem{lemma}[proposition]{Lemma}
\newtheorem{remark}[proposition]{Remark}

\newcommand{\nc}{\newcommand}
\nc{\R}{{\mathbb R}}
\nc{\N}{{\mathbb N}}
\nc{\Z}{{\mathbb Z}}

\nc{\BP}{\mathbb{P}}
\nc{\BE}{\mathbb{E}}
\nc{\BQ}{\mathbb{Q}}
\nc{\bN}{{\mathbf N}}
\nc{\BX}{{\mathbb X}}
\nc{\BY}{{\mathbb Y}}
\nc{\bM}{{\mathbf M}}
\nc{\bF}{{\mathbf F}}
\nc{\bG}{{\mathbf G}}
\nc{\bH}{{\mathbf H}}
\nc{\bW}{{\mathbf W}}
\DeclareMathOperator{\BV}{{\mathbb Var}}
\DeclareMathOperator{\CV}{{\mathbb Cov}}

\nc{\MM}{{\mathcal M}}
\nc{\YY}{{\mathcal Y}}
\nc{\BF}{\mathbf{F}}

\nc{\bea}{\begin{eqnarray}}
\nc{\eea}{\end{eqnarray}}
\nc{\bean}{\begin{eqnarray*}}
\nc{\eean}{\end{eqnarray*}}
\nc{\by}{{\bf y}}
\nc{\bm}{{\bf m}}

\begin{document}

\author{G\"unter Last\footnote{
Institut f\"ur Stochastik, Universit\"at Karlsruhe (TH),
76128 Karlsruhe, Germany. 
Email: last@math.uni-karlsruhe.de}
\ and Mathew D.\ Penrose 
\footnote{
Department of Mathematical Sciences, University of Bath,
Bath BA2 7AY, United Kingdom,
Email: m.d.penrose@bath.ac.uk} 
\footnote{Partially supported by
the Alexander von Humboldt Foundation through
a Friedrich Wilhelm Bessel Research Award.}
}

\title{Poisson process Fock space representation,\\ 
chaos expansion and covariance inequalities}
\date{\today}
\maketitle
\begin{abstract}
\noindent
We consider a Poisson process $\eta$
on an arbitrary measurable space  
with an arbitrary sigma-finite  intensity measure.
We  establish an explicit Fock space representation of square
integrable functions of $\eta$.
As a consequence we identify explicitly, in terms of iterated difference
operators, the integrands in the Wiener-It\^o chaos expansion.
We apply these results to extend well-known variance inequalities
for homogeneous Poisson processes on the line to the general Poisson case.
The Poincar\'e inequality is a special case. Further applications
are covariance identities for Poisson processes on (strictly)
ordered spaces and Harris-FKG-inequalities for monotone functions of $\eta$.
\end{abstract}


\section{Introduction}
\label{secintro}
\setcounter{equation}{0}

The aim of this paper is to develop and to exploit
the basic Fock space structure of a Poisson process $\eta$ on
a measurable space $(\BY,\mathcal{Y})$ with
$\sigma$-finite intensity measure $\lambda$.
In contrast to the literature we do not make any restrictions of
generality, neither imposing a topological structure on
the phase space $(\BY,\mathcal{Y})$, nor assuming the
measure $\lambda$ to be continuous.
Moreover, our results are more explicit than
what was available previously.
We use a probabilistic and non-technical approach
that is based on only a few basic properties
of a Poisson process. 

We now describe the contents of this paper in more detail.
The underlying probability space is denoted by $(\Omega,\mathcal{F},\BP)$.
We interpret the Poisson process $\eta$ as a random element in the space
$\bN:=\bN(\BY)$ of integer-valued $\sigma$-finite measures $\mu$
on $\BY$ equipped with the smallest
$\sigma$-field $\mathcal{N}$
making the mappings $\mu\mapsto\mu(B)$ measurable for all
$B\in\mathcal{Y}$. 
For $y\in\BY$ the {\em difference operator} $D_y$
(also known as the {\em add one cost operator})
is given as follows.
For any measurable $f:\bN\rightarrow\R$ the function
$D_yf$ on $\bN$ is defined by
\begin{align}\label{addone}
D_{y}f(\mu):=f(\mu+\delta_{y})-f(\mu),\quad \mu\in\bN,
\end{align}
where $\delta_{y}$ is the Dirac measure located at a point $y\in\BY$.
 Iterating
this definition,
for $n\ge 2$ and $(y_1,\ldots,y_n)\in\BY^n$
we define a function
$D^{n}_{y_1,\ldots,y_n}f:\bN(\BY)\rightarrow\R$
inductively by
\begin{align}\label{differn}
D^{n}_{y_1,\ldots,y_{n}}f:=D^1_{y_{1}}D^{n-1}_{y_2,\ldots,y_{n}}f,
\end{align}
where $D^1:=D$ and $D^0f = f$. As we shall see,
the operator $D^{n}_{y_1,\ldots,y_n}$ is symmetric in $y_1,\ldots,y_n$.
We define
symmetric and (as it turns out) measurable functions $T_nf$ on $\BY^n$ by
\begin{align}\label{fn2}
T_n f (y_1,\ldots,y_n) :=\BE D^n_{y_1,\ldots,y_n} f(\eta),
\end{align}
and we set $T_0f:=\BE f(\eta)$, whenever these expectations are defined.

By $\langle\cdot,\cdot\rangle_n$ we denote the
scalar product in $L^2(\lambda^n)$ and by
$\|\cdot\|_n$ the associated norm.
Let $\BP_\eta$ denote the distribution of $\eta$.
Then $L^2(\BP_\eta)$ is the space of all measurable 
$f:\bN\rightarrow\R$ satisfying $\BE f(\eta)^2<\infty$.
For $n \in \N$ let $\bH_n$ be the space of symmetric functions in
$L^2(\lambda^n)$, and let $\bH_0 := \R$.
Our first result says that
the mapping $f  \mapsto (T_n(f))_{n\geq 0}$
is an isometry from $L^2(\BP_\eta)$ to the {\em Fock space} 
given by the direct sum of the spaces
$\bH_n$, $n\ge 0$, with $L^2$ norms scaled by $n!^{-1/2}$, as
we describe in more detail in Section \ref{fock}.

\begin{theorem}\label{tcov} Let $f\in L^2(\BP_\eta)$. Then 
\begin{align}\label{varcha}
\BE f(\eta)^2=(\BE f(\eta))^2+
\sum^\infty_{n=1}\frac{1}{n!}\|T_nf\|^2_n.
\end{align}
If also $g\in L^2(\BP_\eta)$, then more generally,
\begin{align}\label{covcha}
\BE f(\eta)g(\eta)=(\BE f(\eta))(\BE g(\eta))+\sum^\infty_{n=1}\frac{1}{n!}
\langle T_nf,T_ng \rangle_n.
\end{align}
\end{theorem}

We shall use Theorem \ref{tcov}
to provide a new proof of the following family of 
inequalities for the variance of Poisson
functionals, which were previously 
given by Houdr\'e and Perez-Abreu \cite{HoPer95}
(but with proof only for a Gaussian analogue)
and by Privault \cite{Priv01} for the case
of normal martingales (including the homogeneous Poisson process on the line).
These estimates involve alternating sums with similar terms  
to those in  \eqref{varcha}, except that we take the expectation
outside the inner product.
Accordingly, let $D^nf(\eta)$ denote the mapping
$(y_1,\ldots,y_n)\mapsto  D^n_{y_1,\ldots,y_n}f(\eta)$.

\begin{theorem}\label{tcovineq} Let $f\in L^2(\BP_\eta)$ and $k\in\N$ be 
such that
\begin{align}\label{Dint}
\BE \|D^n f(\eta)\|_n^2<\infty, \quad n=1,\ldots,2k.
\end{align}
Then 
\begin{align}\label{covineq}
\sum^{2k}_{n=1}\frac{(-1)^{n+1}}{n!}\BE \|D^nf(\eta)\|^2_n
\le \BV[f(\eta)]\le
\sum^{2k-1}_{n=1}\frac{(-1)^{n+1}}{n!}\BE \|D^nf(\eta)\|^2_n.
\end{align}
The first inequality of \eqref{covineq} is an equality if and only if
$D^{2k+1}_{y_1,\ldots,y_{2k+1}}f(\eta) =0$ almost surely,
for $\lambda^{2k+1}$-almost all $(y_1,\ldots,y_{2k+1}) \in \BY^{2k+1}$.
The second inequality of \eqref{covineq} is an equality if and only if
$D^{2k}_{y_1,\ldots,y_{2k}}f(\eta) =0$ almost surely,
for $\lambda^{2k}$-almost all $(y_1,\ldots,y_{2k}) \in \BY^{2k}$.
\end{theorem}
The case $k=1$ of the right hand inequality in \eqref{covineq}
says that
\bea
\BV f( \eta) \leq \BE \int ( f(\eta+ \delta_y) -f(\eta) )^2 \lambda(dy),
\label{poin1}
\eea
and is known as the
{\em Poincar\'e inequality}
for the variance of Poisson functionals.
In the present generality \eqref{poin1} was previously derived by
\cite{Wu00} using an explicit martingale representation.
The same method was used earlier in \cite{Chen85}
to establish the result for infinitely divisible
random vectors with independent components. 
Again our proof is different.

In our opinion the Poincar\'e inequality is a fundamental
property of Poisson processes; for an example of its
application see \cite{HevRei}.
It is related to the well-known {\em Efron-Stein inequality} 
\cite{ES} for  the variance of a symmetric function of $n$
independent $\BY$-valued random variables.

In Section \ref{schaos} we consider for $n\in\N$ the 
{\em multiple Wiener-It\^o integral} $I_n(g)$ of 
a symmetric function $g\in L^2(\lambda^n)$ with respect to
the compensated Poisson process $\hat\eta$
(see \cite{Wiener38,Ito56,NuViv90}).
In fact we shall follow Liebscher \cite{Lieb94} in
defining this integral in the general case, that
is without assuming that $\lambda$ is continuous as in \cite{Ito56,NuViv90}.
For $c\in\R$ we set $I_0(c):=c$. It\^o's \cite{Ito56} 
and Wiener's \cite{Wiener38} famous chaos
expansion of square integrable random variables says that
every function $f \in L^2(\BP_\eta)$ can be decomposed
uniquely as a sum  of variables of the form $I_n(f_n)$
with $(f_n)_{n \geq 0}$ in our Fock space.
The following result identifies each function 
$f_n$ as $1/n!$ times the image of $f$ under the
 mean iterated difference operator $T_n$.

\begin{theorem}\label{tchaos0} 
Let $f\in L^2(\BP_\eta)$. 
Then $T_nf\in L^2(\lambda^n)$, $n\in\N$, and 
\begin{align}\label{chaos}
f(\eta)=\sum^\infty_{n=0}\frac{1}{n!}I_n(T_nf),
\end{align}
where the series converges in $L^2(\BP)$. 
Moreover, if $ g_n \in \bH_n$ for
$n \in \N_0$ satisfy $
f(\eta)=\sum^\infty_{n=0}\frac{1}{n!}I_n(g_n)$ with
convergence in $L^2(\BP)$, then $g_0= \BE f(\eta)$ and
$g_n=T_nf$, $\lambda^n$-almost everywhere on $\BY^n$,
for all $n\in\N$.
\end{theorem}

In the case $\BY=\R$,  Theorem \ref{tchaos0} has been
obtained by Y.\ Ito (\cite{Ito88}, see eqn (7.5) there).
A less explicit version for  the special case of L\'evy processes 
(without Gaussian component)
can be found as Theorem 4 in \cite{Lo05}.
(The proof in \cite{Lo05} does not seem to
justify its usage of iterated stochastic integrals,
see e.g.\ Theorem 18.13 in \cite{Kallenberg} for the Brownian case.)
Our proof is different from these,
and applies to arbitrary $\sigma$-finite intensity measures.
Theorem \ref{tchaos0} and the isometry properties
of stochastic integrals (see \eqref{orth}) show
that the isometry $f\mapsto (T_n(f))_{n\geq 0}$
is in fact a bijection from $L^2(\BP_\eta) $ onto the Fock space.
They could also be used to deduce 
Theorem \ref{tcov}, but we shall proceed in the other
direction, starting with Theorem \ref{tcov} which is
more fundamental. Neither its
formulation nor its proof requires stochastic integration. 
For finite Poisson processes the operators $T_n$  
had been previously used in \cite{MZ960,MolZu00} to
approximate the expectation of a Poisson functional,
while \cite{Bl95} used (in a similar context) a 
closely related operator for more general point processes on the line.

Additional results in Section \ref{schaos} are concerned
with certain derivative and integral operators which are important
in Malliavin calculus on Poisson spaces (see for example \cite{PSTU09}). 
In Theorem \ref{tdiffop}, we provide a
generalization  to arbitrary $\sigma$-finite intensity measure
(and new proof) of a result in \cite{Ito88}
(see also \cite{NuViv90}) identifying the difference operator with 
a stochastic linear derivative operator from 
$L^2(\BP_\eta)$ to $L^2(\BP_\eta \otimes \lambda)$. 
We also consider the
stochastic Kabanov-Skorohod integral 
\cite{Hitsuda72,Kab75,KabSk75} which is a linear operator from
$L^2(\BP_\eta \times \lambda)$ to 
$L^2(\BP_\eta)$ that is dual to the derivative operator.
In Theorem \ref{tequal} we provide a new proof 
of this duality for our more general setting,
and a pathwise interpretation of the Kabanov-Skorohod integral,
using a classical Campbell-type formula due to Mecke \cite{Mecke}
for general Poisson processes.

In Section \ref{scov} we prove some results on covariances,
including the following {\em Harris-FKG inequality}; see also 
\cite{MR96}, \cite{BR06} (treating special Poisson processes)  
\cite{Wu00}, and \cite{GeK97} for more general point processes.
Given $B   \in \YY$,
a function $f:\bN(\BY)\rightarrow\R$ is {\em increasing on} $B$
if $f(\mu+\delta_y)\ge f(\mu)$ for all $\mu\in \bN(\BY)$
and all $y\in B$. It is  
{\em decreasing on} $B$ if $(-f)$ is increasing on $B$.
The formulation we give allows for functions which
are increasing on some parts of $\BY$ and decreasing
on others, which has occasionally been useful; see
Lemma 14 on page 278 of \cite{BR06}, or page 878 of \cite{PS05}.

\begin{theorem}\label{c2}
Suppose $B \in \YY$.
Let $f,g\in L^2(\BP_\eta)$ be increasing on $B$ and decreasing on
$\BY \setminus B$. Then
\begin{align}\label{FKG}
\BE [f(\eta)g(\eta)]\ge (\BE f(\eta))(\BE g(\eta)).
\end{align}
\end{theorem}

We shall derive Theorem \ref{c2} from the following result.
For this result only, we make the extra assumption that
$\BY$ is equipped with a transitive binary relation $<$
such that (i) $\{(y,z):y < z\}$ is a measurable subset of
$\BY^2$ and (ii)
 for any $y,z \in \BY$ at most one
of the relations $y <z $ and $z < y$ can be satisfied,
and (iii) $<$ strictly orders the points of $\BY$ $\lambda$-a.e.,
that is 
\begin{align}\label{diffuse}
\lambda (\BY \setminus \{ z \in \BY: z < y ~~ {\rm or} ~~ y < z \}) 
=0, \quad  y \in \BY.
\end{align}
For any $\mu \in \bN$ let $\mu_y$ denote the restriction of $\mu$
to $y_\downarrow := \{z \in \BY: z < y\}$.
Our final assumption on $<$ is that (iv) $(\mu,y) \mapsto \mu_y$
is a measurable mapping from $\bN \times \BY$ to $\bN$.
In Section \ref{scov} we shall use Theorem \ref{tcov}
to derive the following identity for the 
covariance between two functions of $\eta$.
The theorem requires a version of the
conditional expectation $\BE[D_{y}f(\eta)|\eta_{y}]$
that is jointly measurable in all arguments.
Thanks to the independence properties of a Poisson process
we can and will work with
\begin{align}
\BE[D_{y}f(\eta)|\eta_{y}]:=\int D_{y}f(\eta_{y}+\mu)\Pi^y(d\mu),
\label{956}
\end{align}
where $\Pi^y$ is the distribution of the restriction
of $\eta$ to  $\BY \setminus y_\downarrow$.
By assumption (iv) and Fubini's theorem it
follows that $\Pi^y(\cdot)$ is a kernel, that
is $y\mapsto \Pi^y(A)$ is measurable for all
measurable $A\subset \bN$.

\begin{theorem}\label{tkey} Assume that 
$\BY$ is equipped with a transitive binary relation $<$
satisfying conditions (i)--(iv) above.
For any $f\in L^2(\BP_\eta)$,
\begin{align}\label{sqint}
\BE\int \BE[D_{y}f(\eta)|\eta_{y}]^2\lambda(dy)<\infty,
\end{align}
 and for any $f,g\in L^2(\BP_\eta)$,
\begin{align}\label{Cov}
\CV[f(\eta),g(\eta)]= 
\BE\int \BE[D_{y}f(\eta)|\eta_{y}]\BE[D_{y}g(\eta)|\eta_{y}]\lambda(dy).
\end{align}
\end{theorem}

When $\lambda$ is the product of Lebesgue measure
on $[0,1]$ and a $\sigma$-finite measure on some space $\BX$,
formula \eqref{Cov} is proved in \cite{Chen85}
for special functions $f$, mentioned in \cite{Wu00}
(under an additional assumption on $f$)
and derived (in case of a finite and absolutely continuous 
intensity measure) in \cite{PW08}.
For normal martingales the result is stated in \cite{HoP02}. 
The version for infinitely divisible random vectors can
be found in \cite{Chen85} and \cite{Ho02}.

In Section \ref{sinfdiv} we shall discuss the Poincar\'e and
Harris-FKG inequalities for infinitely divisible random measures.
In Section \ref{sfinite} we describe some of the implications
of our results in the case where $\BY$ is a finite set.

As has been discussed, some of our results already appear in
the literature for special cases such as when the
intensity measure is Lebesgue measure on $\R_+$.
It may be possible to extend  these existing results to
the case where $\BY$ is a Borel space
(a space that is Borel isomorphic to a Borel subset of $[0,1]$),
by considering $\lambda$ as the image of Lebesgue measure
under  a measurable mapping from $\R_+$ to $\BY$ (see Lemma 3.22 of
\cite{Kallenberg}) and making an appropriate change of variables
in the integrals. Nevertheless, we think our direct approach is
worthwhile; as well as being applicable to an arbitrary intensity
measure without any Borel condition, it provides a natural 
approach which relies only on
basic properties of Poisson processes, avoiding the technicalities
of stochastic calculus seen in previous work.

\section{Fock space representation}\label{fock}
\setcounter{equation}{0}

For any $n\in\N$ let $\bH_n$ denote the space
of all measurable functions 
$g:\BY^n\rightarrow\R$ that are square-integrable
with respect to $\lambda^n$ and symmetric $\lambda_n$-a.e.,
equipped with the $L^2(\lambda^n)$ inner product
$\langle \cdot , \cdot \rangle _n$ 
and corresponding norm
 norm $\|\cdot\|_n$ as in Section \ref{secintro}.
Define $\bH_0:=\R$.
Consider the vector space $\bH$ of all
sequences $f=(f_n)_{n\ge 0}$ satisfying
$f_n\in\bH_n$, $n\ge 0$, and
\begin{align}\label{Fock}
\sum_{n=0}^\infty\frac{1}{n!}\|f_n\|_n^2<\infty
\end{align}
where $\|f_0\|_0:=|f_0|$. Equipped with the scalar
product
\begin{align}\label{Fock2}
\langle f,g\rangle_\bH:=\sum_{n=0}^\infty\frac{1}{n!}\langle f_n,g_n\rangle_n,
\quad f=(f_n),g=(g_n)\in\bH, 
\end{align}
$\bH$ becomes a Hilbert space; this is the Fock space that we consider
in this paper. Meyer \cite{Meyer95} gives an introduction
into stochastic calculus on these spaces. 

For any $f\in L^2(\BP_\eta)$ we define
$Tf:=(T_nf)_{n \geq 0}$, where $T_nf$ is given at \eqref{fn2}.
Theorem \ref{tcov} asserts that $Tf\in\bH$ for $f\in L^2(\BP_\eta)$
and
\begin{align}\label{010}
\BE f(\eta)g(\eta)=\langle Tf,Tg \rangle_\bH,\quad f,g\in L^2(\BP_\eta).
\end{align}
We prove these assertions in stages. Note first that
\begin{align}\label{Dsymmetric}
D^n_{y_1,\ldots,y_n}f(\mu)=
\sum_{J \subset \{1,2,\ldots,n\}}(-1)^{n-|J|}
f \Big(\mu+\sum_{j\in J}\delta_{y_j}\Big),
\end{align}
where 
$|J|$ denotes the number of elements of $J$. This
shows that the operator $D^n_{y_1,\ldots,y_n}$
is symmetric in $y_1,\ldots,y_n$, and that
$(\mu,y_1,\ldots,y_n)\mapsto D^n_{y_1,\ldots,y_n}f(\mu)$
is measurable whenever $f:\bN\rightarrow\R$ is measurable.

Let $\mathbf{F}^+$ denote the space of all bounded and measurable
functions $v:\BY\rightarrow\R_+$.
By Lemma 12.2 in \cite{Kallenberg}
the {\em Laplace functional} of $\eta$ is given by
\begin{align}\label{Laplace}
\BE\exp[-\eta(v)]=\exp[-\lambda(1-e^{-v})],\quad v\in \mathbf{F}^+,
\end{align}
where $\mu(v):= \int v d\mu$ for any measure $\mu$ on $\BY$.

Let $\mathcal{Y}_{0}$ be the system of all measurable 
$B\in\mathcal{Y}$ having $\lambda(B)<\infty$.
Let $\mathbf{F}^+_{B}$ denote the space of all those
$v\in \mathbf{F}^+$ vanishing outside $B\in\mathcal{Y}$. 
Let $\mathbf{F}^+_0$ be the space of all functions $v$
that belong to $\mathbf{F}^+_{B}$ for some $B\in\mathcal{Y}_{0}$.
Let $\mathbf{G}$ denote the space of all
(bounded and measurable) functions $g:\bN\rightarrow \R$
of the form
\begin{align}\label{dense4}
g(\mu)=a_1e^{-\mu(v_1)}+\ldots+a_ne^{-\mu(v_n)},
\end{align}
where $n\in\N$, $a_1,\ldots,a_n\in\R$ and $v_1,\ldots,v_n\in \mathbf{F}^+_0$.

\begin{lemma}\label{l9.1}
Relation \eqref{010} holds for $f,g\in\mathbf{G}$.
\end{lemma}
{\sc Proof:} By linearity it suffices to consider
functions $f$ and $g$ of the form
\begin{align*}
f(\mu)=\exp[-\mu(v)],\quad g(\mu)=\exp[-\mu(w)]
\end{align*}
for $v,w\in \mathbf{F}^+_0$. Then we have for $n\ge 1$
that
$$
D^nf(\mu)=\exp[-\mu(v)](e^{-v}-1)^{\otimes n},
$$
where 
$(e^{-v}-1)^{\otimes n}(y_1,\ldots,y_n):=\prod_{i=1}^n (e^{-v(y_i)}-1)$.
From \eqref{Laplace} we obtain that
\begin{align}\label{011}
T_nf=\exp[-\lambda(1-e^{-v})](e^{-v}-1)^{\otimes n}.
\end{align}
Since $v\in\mathbf{F}^+_0$ it follows that
$T_nf\in \bH_n$, $n\ge 0$.
Using \eqref{Laplace}, we obtain that
\begin{align}\label{Lap1}
\BE f(\eta)g(\eta)=\exp[-\lambda(1-e^{-(v+w)})].
\end{align}
On the other hand we have from \eqref{011} 
(putting $\lambda^0(1):=1$) that
\begin{align*}
\sum_{n=0}^\infty\frac{1}{n!}\langle &T_nf,T_ng\rangle_n\\
&=\exp[-\lambda(1-e^{-v})]\exp[-\lambda(1-e^{-w})]
\sum_{n=0}^\infty\frac{1}{n!}\lambda^n(((e^{-v}-1)(e^{-w}-1))^{\otimes n})\\
&=\exp[-\lambda(2-e^{-v}-e^{-w})]
\exp[\lambda((e^{-v}-1)(e^{-w}-1))].
\end{align*}
This equals the right-hand side of \eqref{Lap1}.
\qed

\vspace{0.3cm}
To extend 
\eqref{010} to general 
$f,g\in L^2(\BP_\eta)$ we need two lemmas.

\begin{lemma}\label{lemmdense}\rm 
The set $\bG$ is dense in $L^2(\BP_\eta)$.
\end{lemma}
{\sc Proof:} Let $\mathbf{W}$ be the space of all
bounded measurable $g:\bN\rightarrow\R$ that
can be approximated in $L^2(\BP_\eta)$ by functions
in $\bG$. This space is closed under monotone and
uniformly bounded convergence and contains
the constant functions. 
The space $\mathbf{G}$ is stable under multiplication
and we denote by $\mathcal{N}'$ 
the smallest $\sigma$-field on $\bN$ such that $\mu\mapsto h(\mu)$
is measurable for all $h\in\mathbf{G}$.
A well-known functional version
of the monotone class theorem (see e.g.\ Theorem I.21 in \cite{DeMeyer78})
implies that
$\mathbf{W}$ contains any bounded $\mathcal{N}'$-measurable $g$.
On the other hand we have that
$$
\mu(C)=\lim_{t\to 0+}t^{-1}(1 - e^{-t\mu(C)}),\quad \mu\in \bN,
$$
for any  $C\in \YY$.
Hence $\mu\mapsto\mu(C)$ is $\mathcal{N}'$-measurable 
whenever $C\in\mathcal{Y}_0$.
Since $\lambda$ is $\sigma$-finite, 
for any $C\in\mathcal{Y}$
there is a monotone sequence
$C_k\in\mathcal{Y}_0$, $k\in\N$, with union $C$, so that
 $\mu\mapsto\mu(C)$ is $\mathcal{N}'$-measurable.
Hence $\mathcal{N}'=\mathcal{N}$ and it follows that $\mathbf{W}$ contains all
bounded measurable functions. But then $\mathbf{W}$
is clearly dense in $L^2(\BP_\eta)$ and the proof of
the lemma is complete.\qed

\vspace{0.3cm}
In proving the next lemma, and repeatedly later,
we need to consider {\em factorial moment measures}.
For  $\mu \in \bN(\BY)$ and for $m\in\N$, define the measure
$\mu^{(m)}$ on $\BY^m$  by
\begin{align}\label{mum}
\mu^{(m)}(B):=\idotsint\I_B(y_1,\ldots,y_m)
\Big(\mu - \sum_{j=1}^{m-1} \delta_{y_j}\Big)(dy_m)
\Big(\mu - \sum_{j=1}^{m-2} \delta_{y_j}\Big)(dy_{m-1}) \nonumber \\
\ldots (\mu - \delta_{y_1})(dy_2)\mu(dy_1).
\end{align}
If we can write $\mu = \sum_{i \geq 1} \delta_{x_i}$,
 then $\mu^{(m)}(B)$ counts
the number of $n$-tuples of distinct indices $(i_1,\ldots,i_m)$ 
such that  $(x_{i_1},\ldots,x_{i_k})$ is in $B$, and
$\BE \eta^{(m)} (\cdot)$ is known as the $m$th factorial moment measure
of the Poisson process $\eta$.
A standard tool in the analysis of Poisson driven
stochastic systems is the
formula 
\begin{align}\label{Meckem}\notag
\BE\int h(\eta,y_1,\ldots,y_m)&\eta^{(m)}(d(y_1,\ldots,y_m))\\
&=\BE\int h(\eta+\delta_{y_1}+\ldots+\delta_{y_m},y_1,\ldots,y_m)
\lambda^m(d(y_1,\ldots,y_m)),
\end{align}
for all $h:\bN\times\BY^m\rightarrow[-\infty,\infty]$
for which one (and then also the other) side makes sense. 
When $m=1$, \eqref{Meckem} simplifies
to the following classical formula by Mecke \cite{Mecke}:
\begin{align}\label{Mecke}
\BE\int h(\eta,y)\eta(dy)=\BE\int h(\eta+\delta_{y},y)\lambda(dy).
\end{align}
In the special case where  $\lambda$ is a finite, absolutely continuous
measure on $\R^d$, a proof of \eqref{Meckem} is given
in e.g.\ Theorem 1.6 of \cite{Pen03}, and the argument there can be
extended to the general case. Here, and again later on,
we use the fact that for arbitrary $\sigma$-finite $\lambda$,
the standard proof of existence of a Poisson process $\eta$
with $\sigma$-finite intensity measure
$\lambda$ (see e.g. Theorem 12.7 of \cite{Kallenberg})
shows that there
is a version of this Poisson process process taking the 
form $\sum_{i=1}^{\eta(\BY)} \delta_{X_i}$ for a
sequence of $\BY$-valued random variables $X_i$.

\begin{lemma}\label{lemsubs}\rm Suppose that
$f,f^1,f^2,\ldots\in L^2(\BP_\eta)$ satisfy
$f^k\to f$ in $L^2(\BP_\eta)$ as $k\to\infty$, and that
$h: \bN \to [0,1]$ is measurable. 
Let $n\in\N$, let $C\in\mathcal{Y}_0$ and set $B:=C^n$. 
Then 
\begin{align}\label{013}
\lim_{k\to\infty}\int_B
\BE [| D^n_{y_1,\ldots,y_n} f(\eta) - D^n_{y_1,\ldots,y_n} f^k(\eta)|
h(\eta)  ]  
\lambda^n
(d(y_1,\ldots,y_n))=0.
\end{align}
\end{lemma}
{\sc Proof:} 
By \eqref{Dsymmetric}, the relation \eqref{013} is
implied by the convergence
\begin{align}\label{014}
\lim_{n\to\infty}\int_B
\BE \Big[ \Big|f\Big(\eta+ \sum_{i=1}^m \delta_{y_i} \Big)-
f^k\Big(\eta+ \sum_{i=1}^m \delta_{y_i} \Big) \Big|
h(\eta) \Big]\lambda^n(d(y_1,\ldots,y_n))=0
\end{align}
for all $m\in\{0,\ldots,n\}$. For $m=0$ this is obvious.
Assume $m\in\{1,\ldots,n\}$. Then the integral in \eqref{014}
equals
\begin{align*}
\lambda(C)^{n-m}&\BE \int_{C^m} 
\Big|f\Big(\eta+\sum_{i=1}^m \delta_{y_i} \Big)-
f^k\Big(\eta+ \sum_{i=1}^m \delta_{y_i} \Big) \Big|h(\eta)
\lambda^m(d(y_1,\ldots,y_m))\\
&=\lambda(C)^{n-m}\BE \int_{C^m}  |f(\eta)-f^k(\eta)|
h \Big(\eta - \sum_{i=1}^n \delta_ {y_i} \Big)
\eta^{(m)}(d(y_1, \ldots , y_m)),
\\
& \leq\lambda(C)^{n-m}\BE  |f(\eta)-f^k(\eta)|\eta^{(m)}(C^m),
\end{align*}
where we have used \eqref{Meckem} to get the equality.
By the Cauchy-Schwarz inequality the last expression is bounded
above by 
\begin{align*}
\lambda(C)^{n-m}(\BE (f(\eta)-f^k(\eta))^2)^{1/2}
(\BE[(\eta^{(m)}(C^m))^2])^{1/2}.
\end{align*}
Since any Poisson variable has moments of all orders,
we obtain \eqref{014} and hence the lemma.\qed

\vspace{0.3cm}
Later we will also need the following direct
consequence of \eqref{Mecke}. 

\begin{lemma}\label{l98} Let $f,g:\bN\rightarrow\R$
be measurable functions that coincide $\BP_\eta$-a.e.
Then $D_yf(\mu)=D_yg(\mu)$ for $\BP_\eta\otimes\lambda$-a.e.\
$(\mu,y)$.
\end{lemma}

\vspace{0.3cm}
\noindent
{\sc Proof of Theorem \ref{tcov}:} By linearity and
the polarization identity
$$
4\langle f,g\rangle_\bH=
\langle f+g,f+g\rangle_\bH-\langle f-g,f-g\rangle_\bH
$$
it suffices to prove \eqref{010} for
$f = g \in L^2(\BP)$.
By Lemma \ref{lemmdense}
there are $f^k\in \bG$, $k\in\N$, satisfying
$f^k\to f$ in $L^2(\BP_\eta)$ as $k\to\infty$.
By Lemma \ref{l9.1}, $Tf^k$, $k\in\N$, is a Cauchy sequence
in $\bH$. Let $\tilde f=(\tilde f_n)\in \bH$ be the limit, that is
\begin{align}\label{0lim}
\lim_{k\to\infty}
\sum_{n=0}^\infty\frac{1}{n!}\|T_nf^k-\tilde f_n\|^2_n=0.
\end{align}
Taking the limit in the identity
$\BE f^k(\eta)^2=\langle Tf^k,Tf^k\rangle_\bH$
yields $\BE f(\eta)^2=\langle \tilde f,\tilde f\rangle_\bH$.
Equation \eqref{0lim} implies that $\tilde f_0=\BE f(\eta)=T_0f$.
It remains to show that for any $n\ge 1$,
\begin{align}\label{identify}
\tilde f_n=T_nf\quad \lambda^n\text{-a.e.}
\end{align}
Let $C\in\mathcal{Y}_0$ and $B:=C^n$. 
Let $\lambda^n|_B$
be the restriction of the measure $\lambda^n$ to $B$. 
By \eqref{0lim} 
$T_nf^k$ converges in $L^2(\lambda^n|_B)$ (and hence in
$L^1(\lambda^n|_B)$) to
$\tilde f_n$,
while by the definition \eqref{fn2} of $T_n$, and the case 
$h \equiv 1$ of \eqref{014}, 
$T_nf^k$ converges in $L^1(\lambda^n|_B)$ to 
$T_nf$. Hence 
these $L^1$ limits must be the same almost everywhere, so that 
$\tilde f_n=T_nf$  $\lambda^n$-a.e.\ on $B$. Since $\lambda$
is assumed $\sigma$-finite, this implies \eqref{identify}
and hence the theorem. \qed

\section{Chaos expansion}\label{schaos}
\setcounter{equation}{0}

Given functions $g_i: \BY \to \R$ for $ i = 1,2,\ldots, n$,
define the tensor product function
$\otimes_{i=1}^n g_i$  to be the function from $\BY^n$ to $\R$
which maps each $(y_1,\ldots,y_n)$ to $\prod_{i=1}^n g_i(y_i)$. 
When the functions $g_1,\ldots,g_n$ are all the same
function $g$, we write $g^{\otimes n}$ for this tensor
product function. This is consistent with notation used in
the proof of Lemma \ref{l9.1}.
 
For $n\in\N_0$ and 
$g\in L^2(\lambda^n)$ 
we define
the {\em multiple Wiener-It\^o integral} $I_n(g)$ of $g$ with respect to
the compensated Poisson process $\hat\eta:=\eta-\lambda$
 as follows.
Set $I_0(c):=c$ for $c\in\R$. When $n \geq 1$,
consider first the case where $g$ is {\em of product form},
by which we mean $g= \otimes_{i=1}^n g_i$,
 with each
of the $g_i$ bounded and vanishing outside $B$ for some
$B \in \YY_0$.
For such $g$, set
\begin{align}\label{prodIWI}
I_n(g):=\sum_{J \subset [n]}
(-1)^{n-|J|}\eta^{(|J|)} (\otimes_{j \in J}g_j)
\lambda^{n-|J|}( \otimes_{j' \in [n] \setminus J}g_{j'}),
\end{align}
where
$[n]: =\{1,\ldots,n\}$, and  $|J|$ denotes the number of elements of $J$, 
and $\eta^{(m)}(\cdot)$ denotes integration with respect
to the measure $\eta^{(m)}$  defined by \eqref{mum}, while
$\eta^{(0)} (\otimes_{j \in  \emptyset}g_j) :=1$ and
$\lambda^0 (\otimes_{j \in \emptyset} g_j) := 1$.
In the special case where the functions $g_i$ are all the same
function $h$, the formula \eqref{prodIWI} simplifies to
\begin{align}\label{symprodIWI}
I_n( h^{\otimes n})=\sum_{k =0}^n\binom{n}{k}(-1)^{n-k}
\eta^{(k)} (h^{\otimes k})(\lambda(h))^{n-k}.
\end{align}

Recall \cite{Ogura72,Ito88,BHJ92} that the Charlier polynomials
$C_n(\lambda;\cdot)$, $n\in\bN_0$, are a family of orthogonal
polynomials for
the Poisson distribution with parameter $\lambda\ge 0$ defined by 
$$
C_{n}(\lambda;x) = \sum_{k=0}^n \binom{n}{k} (-1)^{n-k} \lambda^{-k} (x)_{k}
$$
where $(x)_j$ is the descending factorial $x(x-1)\cdot\ldots\cdot(x-j+1)$
with $(x)_0$ interpreted as 1. Assume that the function $g$ in
\eqref{prodIWI} is of the form 
$g=\I_{B_1}^{\otimes{m_1}}\otimes\ldots\otimes \I_{B_k}^{\otimes{m_k}}$,
where $B_1,\ldots,B_k\in\mathcal{Y}_0$ are pairwise disjoint
and $m_1+\ldots+m_k=n$. Then 
\begin{align*}
I_n(g)=\sum_{j_1=0}^{m_1} \ldots \sum_{j_k=0}^{m_k}
\prod_{i=1}^k 
\binom{m_i}{j_i}(\eta(B_i))_{j_i}(-\lambda(B_i))^{m_i-j_i}.
\end{align*}
This can be written as
\begin{align}\label{stcharl}
I_n\big(\I_{B_1}^{\otimes{m_1}}\otimes\ldots\otimes \I_{B_k}^{\otimes{m_k}}\big)
=\prod_{i=1}^k\lambda(B_i)^{m_i}C_{m_i}(\lambda(B_i);\eta(B_i)),
\end{align}
at least if $\eta(B_i)<\infty$ for $i=1,\ldots,k$, an event
with probability $1$. Ogura \cite{Ogura72}
used this formula to define the Wiener-It\^o integral 
for a homogeneous Poisson process on the line.
Liebscher \cite{Lieb94} generalized this approach
to Poisson processes on a complete separable metric space
with locally finite intensity measure.

We extend the definition of $I_n(g)$ by linearity to those 
functions $g$ which can be expressed as a finite  sum
of functions of product form; we shall say that such $g$ are of
{\em sum-product form}.
The extension is well defined because each term of the sum in
the right hand side of \eqref{prodIWI} is $\pm 1$ times
the integral of $g$ with respect to   a certain
measure on $\BY^m$, and hence is linear in $g$.

Let $\Sigma_n$ denote the set of all permutations of 
$[n]$, and for $g \in L^2(\lambda^n)$ define the
{\em symmetrization} $\tilde g$ of $g$ by 
\begin{align} \label{symdef}
\tilde g(y_1,\ldots,y_n):=\frac{1}{n!}\sum_{\pi \in \Sigma_n} 
g(y_{\pi(1)},\ldots,y_{\pi(n)}).
\end{align}
By \eqref{prodIWI}, for $g$ of product form
$I_n(g)$ is invariant under permutations
of the functions $g_i$ in the tensor product.
Hence $I_n(g) = I_n(\tilde g)$ for all $g$ of sum-product form.

We shall show in Lemma  \ref{l34} below that
if $g $ and $h $ are functions of sum-product form on  $\BY^n$ and
on $ \BY^m$ respectively, we have the isometry relation
\begin{align}\label{orth}
\BE I_m(g)I_n(h)=\I\{m=n\}m!
\langle \tilde g , \tilde h \rangle_n, \quad m,n\in\N_0.
\end{align}
Since functions of sum-product form are dense in  $L^2(\lambda^n)$, 
we can (and do) extend the definition
of $I_n(g)$ to general
$g \in L^2(\lambda^n)$ by isometry.
It follows from the isometry that
 $I_m (g) = I_m(\tilde g)$ for all $g \in L^2(\lambda^m)$,
and that $I_m(g)$ is linear in $g$, and that
\eqref{orth} holds for all $g \in L^2(\lambda^m)$
and $h\in L^2(\lambda^n)$, since all these properties
were already established for functions of sum-product form.

The  proofs of \eqref{orth} in the literature (see \cite{Ito56}, \cite{Surg84})
assume further structure on the measure space 
$(\BY,\YY,\lambda)$, typically  including
diffuseness or continuity of $\lambda$, and in some cases,
topological assumptions on $\BY$. We are not making such assumptions
here, so we provide a new proof of  \eqref{orth}, as follows.
An alternative proof can be based on \eqref{stcharl}
and the orthogonality properties of the  Charlier polynomials,
see \cite{Lieb94}. 

\begin{lemma}\label{l34} Let $m,n\in\N_0$, and
suppose $g:\BY^m \to \R$ and  $h:\BY^n \to \R$ are of sum-product form.
Then $I_m(g),I_n(h)\in L^2(\BP)$, and \eqref{orth} holds.
\end{lemma} 
{\sc Proof:}
First consider $g$ and $h$  of product form.
Suppose 
$n > 0$ and $m >0$, and suppose
$g = \otimes_{j=1}^n g_j$ and
$h = \otimes_{k=1}^m h_k$,
where all the functions $g_j$ and $h_k$ are bounded and vanish outside of
some $B \in \YY_0$.
Then
\begin{align*}
\BE I_m(g) I_n(h)&=\sum_{J \subset [m]} \sum_{K \subset [n]}
(-1)^{m+ n - |J|-|K|}   
\Big( \prod_{j' \in [m ] \setminus J} \lambda(g_{j'}) \Big)
\Big( \prod_{k' \in [n ] \setminus J} \lambda(g_{k'}) \Big) \\
&\quad \times\BE \eta^{(|J|)} ( \otimes_{j \in J} g_j)
\eta^{(|K|)} ( \otimes_{k \in K} h_k).
\end{align*}
Suppose, in each term of the preceding sum, we were
to replace the expectation by
$
\BE [ \eta^{(|J|+|K|)} ( (\otimes_{j \in J} g_j)
\otimes
 (\otimes_{k \in K} h_k) )]
$. Then the modified sum would come to zero,  because
by \eqref{Meckem}, the modified expectation comes to
$(\prod_{j \in  J} \lambda(\eta_j) ) (  
\prod_{k \in  K} \lambda(\eta_k) )$, so in the modified
sum  each term is of the form $(-1)^{m +n -|J|-|K|}
(\prod_{j=1}^m \lambda  (g_j)) 
(\prod_{k=1}^n \lambda  (h_j))$ (this also shows that \eqref{orth} holds
when one of $m$ and $n$ is 0.) Therefore,
\begin{align}\label{090714b}
\BE I_m(g) I_n(h)&=\sum_{J \subset [m]} \sum_{K \subset [n]}
(-1)^{m+ n - |J|-|K|}   
\Big( \prod_{j' \in [m ] \setminus J} \lambda(g_{j'}) \Big)
\Big( \prod_{k' \in [n ] \setminus K} \lambda(g_{k'}) \Big)
\\ \notag
&\quad\times\BE [ \eta^{(|J|)} ( \otimes_{j \in J} g_j)
\eta^{(|K|)} ( \otimes_{k \in K} h_k) - 
 \eta^{(|J|+|K|)} ( (\otimes_{j \in J} g_j) 
\otimes (
 \otimes_{k \in K} h_k) )].
\end{align}
Suppose we write the restriction of $\eta$ to $B$ as
$\sum_{i=1}^N \delta_{ X_i}$ (for any Poisson process this is possible,
as remarked earlier).  
Then
 $\eta^{(|J|)} (\prod_{j \in J} g_j)$
is the sum, over all $|J|$-tuples of distinct points $X_i$,
$i \in [N]$,
of the product of the values of $g_j$ at those points.
Similarly, $\eta^{(|K|)} (\otimes_{k \in K} h_k)$
is the sum, over all $|K|$-tuples of distinct points $X_i$,
of the product of the values of $h_j$ at those points.
When multiplied together, the $j$-tuple and the $k$-tuple
need to have at least one element in common for this product
to be different from the corresponding term in $\eta^{(|J|+|K|)}(
(\otimes_{j \in J} g_j ) \otimes ( \otimes_{k \in K} h_k))$. 
For example, if $J= K =\{1,2\}$
we have
\begin{align*}
\eta^{(2)}(g_1 \otimes g_2)& 
\eta^{(2)}(  h_1 \otimes h_2) - 
\eta^{(4)}(   g_1 \otimes g_2  \otimes   h_1 \otimes h_2) \\
=&  \eta^{(3)} ( g_1 h_1 \otimes g_2 \otimes h_2 +
 g_1 h_2 \otimes g_2 \otimes h_1 +
 g_2 h_1 \otimes g_1 \otimes h_2 +
 g_2 h_2 \otimes g_1 \otimes h_1 )\\
&+ \eta^{(2)}(   g_1 h_1 \otimes g_2 h_2 +  g_1 h_2 \otimes g_2 h_1) 
\end{align*}
where the product $g_i h_j$ is defined pointwise, that is
$g_ih_j(y) := g_i(y) h_j(y)$ for $y \in \BY$.
In general, for each matching (bijection) $\varphi$  of a nonempty
subset $\{j_1,\ldots,j_\alpha\}$ of $J$  
to  a  subset  
of $K$, writing $ k_i $ for $ \varphi(j_i) $
we get a contribution
to   the expression inside the expectation in the $(J,K)$th term
of \eqref{090714b}
which is of the form
$$
\eta^{(|J|+|K|- \alpha)}
( (  \otimes_{i=1}^\alpha g_{j_i} h_{k_i})
\otimes  (\otimes_{j \in J \setminus \{j_1,\ldots,j_\alpha\}} 
g_j )
\otimes (  \otimes_{k \in K \setminus \{k_1,\ldots,k_\alpha\}}  
h_k))
$$
so that when one takes the expectation using \eqref{Meckem},
and  multiplies by the remaining factors
of $\lambda(g_{j'}) \lambda(h_{k'}) $ appearing in
this term of the right hand side of \eqref{090714b},
one ends up with a contribution of
\begin{align*}
(-1)^{m+ n - |J|-|K|}   
 \Big(  \prod_{i=1}^\alpha \lambda( g_{j_i} h_{k_i}) \Big)
\times  \Big( \prod_{j \in [m] \setminus \{j_1,\ldots,j_\alpha\}} 
 \lambda(g_j) \Big)
\times  \Big( \prod_{k \in [n] \setminus \{k_1,\ldots,k_\alpha\}} 
\lambda(h_k) \Big)
\end{align*}
which depends on $J$ and $K$ only through the sign factor.
Writing 
 $J = \{j_1,\ldots,j_\alpha\}\cup J'$ and 
 $K = \{k_1,\ldots,k_\alpha\}\cup K'$, we have that
 the total contribution to the right hand side of
\eqref{090714b} from a given matching $\varphi$ is given
by
\begin{align*}
&\Big(\prod_{i=1}^\alpha \lambda ( g_{j_i} h_{k_i}) \Big)
\times   \Big( \prod_{j \in [m] \setminus \{j_1,\ldots,j_\alpha\}} 
 \lambda(g_j) \Big)
\times  \Big( \prod_{k \in [n] \setminus \{k_1,\ldots,k_\alpha\}} 
\lambda(h_k) \Big)\\
&\quad\times \sum_{J' \subset [m]  \setminus \{j_1,\ldots,j_\alpha \} }
\sum_{K'   \subset  [n] \setminus \{k_1,\ldots,k_\alpha \} } 
(-1)^{m+ n - |J'|-|K'| - 2 \alpha}   
\end{align*}
and this comes to zero, except in the case where
$\alpha = m = n$.

Hence, all matchings contribute zero to \eqref{090714b} unless
$m=n$, and for this case there are $n!$ matchings
$\varphi$ having $\alpha = n$,  namely
the permutations of $[n]$, so that 
\bean
\BE \left[ I_n(g) I_n(h) \right] = 
\sum_{\varphi \in \Sigma_n} 
\prod_{i=1}^n \lambda (g_{i } h_{\varphi(i)} ) 
=
\sum_{\varphi \in \Sigma_n} 
\prod_{i=1}^n \langle g_{i }, h_{\varphi(i)} \rangle_n. 
\eean
With the symmetrization $\tilde{g}$ defined at \eqref{symdef}, we have
by linearity that
\begin{align*}
\BE  I_n(\tilde{g}) I_n(\tilde{h})=n!^{-2}
\sum_{\pi \in \Sigma_n} \sum_{\sigma \in \Sigma_n} 
\BE [ I_n ( \otimes_{i=1}^n g_{\pi(i)})
 I_n ( \otimes_{j=1}^n h_{\sigma(j)})]
=\sum_{\varphi \in \Sigma_n} \prod_{i=1}^n
\langle g_{i } ,h_{\varphi(i)} \rangle_n 
\end{align*}
whereas 
\begin{align*}
\langle \tilde{g}, \tilde{h} \rangle_n 
= n!^{-2}  
 \sum_{\pi \in \Sigma_n} 
 \sum_{\sigma \in \Sigma_n} 
\langle g_{\pi(i)}, h_{\sigma(j)} \rangle_n
= n!^{-1} \sum_{\varphi \in \Sigma_n} 
\langle g_{i}, h_{\varphi(i)} \rangle_n 
\end{align*}
so that \eqref{orth} holds for this case. We can then
extend by linearity to all $f$ and $g$ of sum-product form.\qed

\vspace{0.5cm}

The proof of Theorem \ref{tchaos0} requires
the following key lemma.

\begin{lemma}\label{l3} Let $f(\mu):=e^{-\mu(v)}$, $\mu\in\bN(\BY)$,
where $v:\BY\rightarrow\R_+$
is a measurable function vanishing outside a set $B\in\mathcal{Y}_0$.
Then \eqref{chaos} holds $\BP$-a.s.\ and in $L^2(\BP)$.
\end{lemma}
{\sc Proof:} 
By \eqref{Laplace} and \eqref{011}
the right-hand side of \eqref{chaos} equals the formal sum
\begin{align}
I:=\exp[-\lambda(1-e^{-v})]
+\exp[-\lambda(1-e^{-v})]\sum^\infty_{n=1}\frac{1}{n!}
I_n((e^{-v}-1)^{\otimes n}).
\label{chaos9}
\end{align}
The function
$(e^{-v}-1)^{\otimes n}$ is of product form.
Using the pathwise definition \eqref{symprodIWI} we obtain that almost surely
\begin{align}
I&=\exp[-\lambda(1-e^{-v})]\sum^\infty_{n=0}\frac{1}{n!}
\sum^n_{k=0}\binom{n}{k}
\eta^{(k)}((e^{-v}-1)^{\otimes k})(\lambda(1-e^{-v}))^{n-k}
\nonumber \\
&=\exp[-\lambda(1-e^{-v})]
\sum^\infty_{k=0}\frac{1}{k!}\eta^{(k)}((e^{-v}-1)^{\otimes k})
\sum^\infty_{n=k}\frac{1}{(n-k)!}(\lambda(1-e^{-v}))^{n-k}
\nonumber \\
&=\sum^{N}_{k=0}\frac{1}{k!}\eta^{(k)}((e^{-v}-1)^{\otimes k}),
\label{0825a}
\end{align}
where $N:=\eta(B)$. Writing
  $\delta_{X_1}+\ldots+\delta_{X_N}$
for the restriction of $\eta$ to $B$,
we have almost surely that
%
%
\begin{align*}
I=
\sum_{J\subset  \{1,\ldots,N\} }
\prod_{i\in J}(e^{-v(X_i)}-1)
&=\prod^N_{i=1}e^{-v(X_i)}=e^{-\eta(v)},
\end{align*}
and hence \eqref{chaos} holds with almost sure convergence of the series.
To demonstrate that convergence also holds in $L^2(\BP)$,
let the partial sum $I(m)$ be given by the right hand side \eqref{chaos9} with  
the series terminated at $n=m$. Then since $\lambda(1-e^{-v})$
is nonnegative and $|1-e^{-v(y)}| \leq 1$ for all $y$,
a similar argument to \eqref{0825a} yields
\begin{align*}
|I(m)| & \leq  
\sum^{\min(N,m)}_{k=0}\frac{1}{k!} | \eta^{(k)}((e^{-v}-1)^{\otimes k}) | \\
 & \leq  
\sum^{N}_{k=0}\frac{N(N-1) \cdots (N-k+1)}{k!} 
= 2^N.
\end{align*}
Since $2^N$ has finite moments of all orders, by dominated
convergence the series
\eqref{chaos9} (and hence \eqref{chaos}) converges in $L^2(\BP)$.
\vspace{0.3cm}
\noindent
{\sc Proof of Theorem \ref{tchaos0}:} Let $f \in L^2(\BP_\eta)$
and define $T_n f$ for $n\in\N_0$ by \eqref{fn2}. 
By \eqref{orth} and  Theorem \ref{tcov}, 
$$
\sum^\infty_{n=0}\BE\Big(\frac{1}{n!}I_n(T_nf)\Big)^2
=\sum^\infty_{n=0}\frac{1}{n!}\|T_nf\|_n^2 = \BE f(\eta)^2 <\infty.
$$
Hence the infinite series of orthogonal terms
$$
X:=\sum^\infty_{n=0}\frac{1}{n!}I_n(T_nf)
$$
converges in $L^2(\BP)$. Let $h \in\bG$,  
where $\bG$ was defined at \eqref{dense4}. 
By Lemma \ref{l3} and linearity of $I_n(\cdot)$
the sum $\sum_{n=0}^\infty \frac{1}{n!} I_n(T_n h)$
converges in $L^2(\BP)$ to $ h(\eta)$.  
Using \eqref{orth} followed by Theorem \ref{tcov} yields
\begin{align*}
\BE (h(\eta) - X)^2=\sum^\infty_{n=0}\frac{1}{n!}\|T_nh - T_n f\|_{n}
= \BE(f(\eta) -h(\eta))^2.
\end{align*}
Hence if $\BE(f(\eta) - h(\eta))^2$ is small, then so is
$\BE(f(\eta) - X)^2$.
Since $\bG$ dense in $L^2(\BP_\eta)$
by  Lemma \ref{lemmdense}, it follows that  $f(\eta) = X$ almost
surely.

To prove the uniqueness, suppose that 
also $ g_n \in \bH_n$ for
$n \in \N_0$ are such that
$\sum^\infty_{n=0}\frac{1}{n!}I_n(g_n)$ 
 converges in $L^2(\BP)$ to $f(\eta)$.
By taking expectations  we must have $g_0=\BE f(\eta)= T_0 f$.
For $n\ge 1$ and $h \in \bH_n$, by \eqref{orth} and \eqref{chaos}
we have
$$
\BE f(\eta) I_n(h) = \BE I_n(T_n f) I_n(h)= n!\langle T_n f, h \rangle_n 
$$ 
and similarly with $T_n f$ replaced by $g_n$, so that
$\langle T_nf -g_n, h \rangle_n =0$. Putting $h= T_nf-g_n$
gives $\| T_nf -g_n \|_n =0$ for each $n$, completing   
the proof of the theorem. \qed
\vspace{0.3cm}

%

\vspace{0.3cm}
We proceed with proving that the pathwise
defined difference operator $D_y$ coincides
with a derivative operator acting on square integrable 
$\sigma(\eta)$-measurable random variables. For $n\in\N$, $f\in\bH_n$,
and $y\in\BY$ we define
$f_y:\BY^{n-1}\rightarrow \R$ by
$f_y(y_1,\ldots,y_{n-1}):=f(y_1,\ldots,y_{n-1},y)$.
Whenever $f_y\in L^2(\lambda^{n-1})$ (which is the case
for $\lambda$-a.e.\ $y$) we define $I_{n-1}f(y):=I_{n-1}(f_y)$.
Otherwise we set $I_{n-1}f(y):=0$. 
We choose a version of
$I_{n-1}f(y)$ that is jointly  measurable in $\omega$ and $y$.
Strictly speaking, we claim that there is a 
$h\in L^2(\BP\otimes\lambda^n)$ and a 
$B\in\mathcal{Y}$ such that $\lambda(\BY\setminus B)=0$ and
$I_{n-1}f(y)=h(\cdot,y)$ $\BP$-a.s.\ for any $y\in B$.
If $f$ was a function of sum-product form
then by using \eqref{prodIWI} one could see directly that $I_{n-1}f$
was jointly  measurable. In general we approximate 
$f$ in $L^2(\lambda^n)$ by symmetric functions $f^k$, $k\in\N$, of
sum-product form. Since by \eqref{orth}
\begin{align*}
\int\BE |I_{n-1}f^k(y)-I_{n-1}f^l(y)|^2\lambda(dy)=
\int\|f^k_y-f^l_y\|_{n-1}^2\lambda(dy)=\|f^k-f^l\|_{n}^2,
\quad k,l\in\N,
\end{align*}
we can take $h$ as the $L^2$-limit of the Cauchy sequence $I_{n-1}f^k(y)$.
Next we can choose $B_1\in\mathcal{Y}$ and a subsequence $J\subset\N$
such that $\lambda(\BY\setminus B_1)=0$  and 
$\BE |I_{n-1}f^k(y)-h(\cdot,y)|^2\to 0$ as $k\to\infty$ along $J$
for all $y\in B_1$. On the other hand we may also choose
$B_2\in\mathcal{Y}$ and a subsequence $J'\subset J$
such that $\lambda(\BY\setminus B_2)=0$  and 
$\|f^k_y-f_y\|^2_{n-1}\to 0$ as $k\to\infty$ along $J'$
for all $y\in B_2$. But then
$\BE |I_{n-1}f^k(y)-I_{n-1}f(y)|^2\to 0$ as $k\to\infty$ along $J'$,
implying that
$I_{n-1}f(y)=h(\cdot,y)$ $\BP$-a.s.\ for any $y\in B_1\cap B_2$.

Given $f\in L^2(\BP_\eta)$, define $f_n = \frac{1}{n!} T_nf \in \bH_n$,
so by Theorem \ref{tchaos0},
\begin{align}\label{chaos2}
f(\eta)=\sum^\infty_{n=0}I_n(f_n) 
\end{align}
is the chaotic expansion of $f$ (with $L^2(\BP)$ convergence).
We then define
\begin{align}\label{D'}
D'_yf (\eta) :=\sum^\infty_{n=1} nI_{n-1}f_n(y),
\end{align}
provided that
\begin{align}\label{conv}
\sum^\infty_{n=1}n\cdot n!\int f^2_nd\lambda^n<\infty.
\end{align}
In this case, by \eqref{orth},
\begin{align*}
\BE\int (D'_yf (\eta) )^2\lambda(dy)
=\int\Big(\sum^\infty_{n=1}n\cdot n!\|(f_n)_y\|^2_{n-1}\Big)\lambda(dy)
=\sum^\infty_{n=1}n\cdot n!\int \|f_n\|^2_n.
\end{align*}
Therefore we can interpret $D'$ as a linear {\em derivative operator}
from $L^2(\BP_\eta)$ to $L^2(\BP_\eta\otimes\lambda)$, see e.g.\ \cite{NuViv90}.
The following result generalizes Theorem 6.5 in \cite{Ito88}
(see also Theorem 6.2 in \cite{NuViv90}).

\begin{theorem}\label{tdiffop} Let $f\in L^2(\BP_\eta)$ have
a chaotic expansion \eqref{chaos2} satisfying \eqref{conv}. Then
\begin{align}\label{11}
D_yf(\eta)=D'_yf( \eta) \quad \BP\text{-a.s.},
\lambda\text{-a.e.\ $y$}.
\end{align}
\end{theorem}

Before proving Theorem \ref{tdiffop}, it is convenient
to introduce the dual operator of $D'$.
Let $h\in L^2(\BP_\eta\otimes\lambda)$. Then $h(\cdot,y)\in L^2(\BP_\eta)$
for $\lambda$-a.e.\ $y$.
Define
\begin{align}
\label{0811a}
\tilde h_n(y_1,\ldots,y_{n+1})=
\frac{1}{(n+1)!}
\sum^{n+1}_{i=1}
\BE D^{n}_{y_1,\ldots,y_{i-1},y_{i+1},\ldots,y_n}f(\eta,y_i) .
\end{align}
From 
Theorem \ref{tcov} we obtain
that $\tilde h_n\in \bH_{n+1}$ and we can define
the {\em Kabanov-Skorohod integral} \cite{Hitsuda72,Kab75,Skor75,KabSk75}
of $h$, denoted
$\delta(h)$, by
\begin{align}\label{Skorint}
\delta(h):=\sum^\infty_{n=0} I_{n+1}(\tilde h_n),
\end{align}
which converges in $L^2(\BP)$
 provided that
\begin{align}\label{domainS}
\sum^\infty_{n=0}(n+1)!\int \tilde h_n^2d\lambda^{n+1}<\infty.
\end{align}

The following duality relation is a special case of
Proposition 4.2 in \cite{NuViv90} applying to general
Fock spaces. We give the short proof for completeness.

\begin{proposition}\label{Spath} Let $f\in L^2(\BP_\eta)$ 
have a chaotic expansion \eqref{chaos2}  satisfying \eqref{conv}
and let $h\in L^2(\BP_\eta\otimes\lambda)$
be such that \eqref{domainS} holds.
Then
\begin{align}\label{adj2}
\BE\int D'_yf(\eta)h(\eta,y)\lambda(dy)=\BE f(\eta)\delta(h).
\end{align}
\end{proposition}
{\sc Proof:} For $y \in \BY$ with $h(\cdot,y) \in L^2(\BP_\eta)$, 
the function $T_n h(\cdot,y) \in \bH_n$ is defined,
and by Theorem \ref{tchaos0} the sum
$\sum^\infty_{n=0} \frac{1}{n!} I_n(T_n h(\cdot,y))$, converges
in $L^2(\BP)$ to
$h(\eta,y)$.
Also,
by the explicit formula \eqref{Dsymmetric},
the function
$(y_1,\ldots,y_n,y)\mapsto  \frac{1}{n!} T_nh(\cdot,y)(y_1,\ldots,y_n)$
is measurable. Moreover $\tilde h_n$, given by \eqref{0811a}, is
the  symmetrization 
of this function, as defined at \eqref{symdef}.
By \eqref{orth},
\begin{align*}
\BE\int D'_yf(\eta)h(\eta,y)\lambda(dy)
=\sum^\infty_{n=1}\int \langle (f_n)_y,T_{n-1}h(\cdot,y)\rangle_{n-1}\lambda(dy)
=\sum^\infty_{n=1}n! \langle f_n,\tilde h_{n-1}\rangle_n,
\end{align*}
where we have used the fact that $f_{n}$ is a symmetric function.
By definition \eqref{Skorint} and \eqref{orth}, 
the last series coincides with $\BE f(\eta)\delta(h)$.\qed


\vspace{0.3cm}
{\sc Proof of Theorem \ref{tdiffop}:} First consider
the case with $f(\mu) = e^{-\mu(v)} $ with 
$v \in \BF_0^+$. Then $n!f_n = T_nf$ is given by \eqref{011}.
Given $n \in \N$,
\begin{align*}
n \cdot n!
\int f_n^2 d \lambda^n  = \frac{1}{(n-1)!} \exp[2 \lambda(e^{-v}-1)] 
(\lambda((e^{-v}-1)^2))^{n}
\end{align*}
which is summable in $n$, so \eqref{conv} holds in this case.
Also, in this case, $D_yf(\eta) = (e^{v(y)}-1) f(\eta)$ by
\eqref{addone}, while
 $
(f_{n})_y = 
(e^{-v(y)} -1) n^{-1} f_{n-1} 
$
 so that by \eqref{D'},
$$
D'_y f( \eta)  = \sum_{n=1}^\infty  
(e^{-v(y)} -1) I_{n-1} (   f_{n-1} ) = 
(e^{-v(y)} -1) f  (\eta) 
$$
where the last inequality is from Lemma  \ref{l3} again.
Thus \eqref{11} holds for $f$ of this form.
By linearity this extends to all elements of $\bG$. 

Let us now consider the general case. 
Choose $g_k\in\bG$, $k\in\N$, such that $g_k\to f$ in $L^2(\BP_\eta)$
as $k\to\infty$, see Lemma \ref{lemmdense}.
Define $h\in L^2(\BP_\eta\otimes\lambda)$
by $h(\mu,y):=h'(\mu)\I_B(y)$, where $h'$ is as in
Lemma \ref{l3} and
 $B\in\mathcal{Y}_0$.  From
 Lemma \ref{l3} it is easy to see that \eqref{domainS}
holds. Therefore we obtain from Proposition \ref{Spath} 
and the linearity of the operator $D'$ that
\begin{align}\label{adj7}
\BE\int (D'_yf (\eta) -D'_yg_k (\eta) )h(\eta,y)\lambda(dy)=\BE (f(\eta)-g_k(\eta))\delta(h)
\to 0\quad \text{as $k\to\infty$}.
\end{align}
On the other hand,
$$
\BE\int (D_yf(\eta)-D_yg_k(\eta))h(\eta,y)\lambda(dy)
= \int_B\BE [
 (D_yf(\eta)-D_yg_k(\eta))h'(\eta)]\lambda(dy),
$$
and  by the case $n=1$ of Lemma \ref{lemsubs},
this tends to zero as $k \to \infty$.
%
Since $D'_yg_k=D_yg_k$ a.s.\ for $\lambda$-a.e.\ $y$
we obtain from \eqref{adj7} that
\begin{align}\label{49}
\BE\int (D'_yf)h(\eta,y)\lambda(dy)=\BE\int (D_yf(\eta))h(\eta,y)\lambda(dy).
\end{align}
By Lemma \ref{lemmdense}, the linear combinations
of the functions $h$ considered above are dense in
$L^2(\BP_\eta \otimes \lambda)$, and by linearity 
\eqref{49} carries through to $h$ in this dense class of functions too, 
so we may conclude that the assertion \eqref{11} holds.\qed

\vspace{0.3cm}
Next we derive a pathwise interpretation of the Kabanov-Skorohod integral. 
For $h\in L^1(\BP_\eta\otimes\lambda)$ we define
\begin{align}\label{pathSk}
\delta'(h):=\int h(\eta - \delta_y,y)\eta(dy)
-\int h(\eta,y)\lambda(dy).
\end{align}
It turns out that the Kabanov-Skorohod integral 
and the operator $\delta'$ coincide
on the intersection of their domains:

\begin{theorem}\label{tequal} Let
$h\in L^1(\BP_\eta\otimes\lambda)\cap L^2(\BP_\eta\otimes\lambda)$
and assume that \eqref{domainS} holds. Then
$\delta(h)=\delta'(h)$ $\BP$-a.s.
\end{theorem}
{\sc Proof:} The Mecke equation \eqref{Mecke} shows that 
$\BE\int |h(\eta - \delta_y,y)|\eta(dy)<\infty$
as well as
\begin{align}\label{adj3}
\BE\int D_yf(\eta)h(\eta,y)\lambda(dy)=\BE f(\eta)\delta'(h),
\end{align}
whenever $f:\bN\rightarrow \R$ is measurable and bounded.
Therefore we obtain from Theorem \ref{tdiffop} and Proposition \ref{Spath}
that $\BE f(\eta)\delta'(h)=\BE f(\eta)\delta(h)$ provided
that $f$ satisfies \eqref{conv}. By Lemma \ref{lemmdense} 
 the space of such
bounded functions is dense in $L^2(\BP_\eta)$, so  we
may conclude that the assertion holds.\qed

\vspace{0.3cm}
The duality relation \eqref{adj3} was observed in \cite{Pic96b},
at least in case $\lambda$ is diffuse.
Therefore Theorem \ref{tequal} is implicit in this work.

\section{Variance inequalities}\label{ScovPoiss}
\setcounter{equation}{0}

In this section we prove Theorem \ref{tcovineq}, and also
give a further 
set of variance inequalities
in Theorem \ref{tpoisecond}. 
To help us identify the cases where these equalities are strict, 
 we first 
give a criterion, in terms of the difference operator, 
for the  chaotic expansion \eqref{chaos2} of a function $f \in L^2(\BP_\eta)$
to  terminate after $k$ steps. 

\begin{proposition}\label{pfinchaos}
Suppose $f \in L^2(\BP_\eta)$  and $k \in \N_0$.
Then $f $ satisfies 
\begin{align}\label{678a}
f(\eta)=I_0(f_0)+\ldots+I_k(f_k) \quad  \BP\text{-a.s.}
\end{align}
for some $f_j \in \bH_j$, $0\leq j \leq k$, if and only if 
$D^{k+1}_{y_1,\ldots,y_{k+1}}f(\eta) =0$ almost surely
for $\lambda^{k+1}$-a.e. $(y_1,\ldots,y_{k+1})$,  
in which case $f_0=\BE f(\eta)$ and $f_j$ and $T_jf$ are equal $\lambda^j$-a.e.\ 
for $j =1,\ldots,k$.
\end{proposition}
{\sc Proof:} 
First suppose $D^{k+1}_{y_1,\ldots,y_{k+1}}f(\eta) =0$  almost surely
for $\lambda^{k+1}$-a.e.\ $(y_1,\ldots,y_{k+1})$.  
Then given such $y_1,\ldots,y_{k+1}$,
for any $m \in \N$ with $m \geq k+2$ 
we have almost surely
$D^{m}_{y_1,\ldots,y_{m}}(f) =0$  
for any  $y_{k+2},\ldots,y_{m}$, and taking expectations
shows that $T^{m}f(y_1,\ldots,y_{m})$
for any $m \geq k+1$ and $\lambda^{m}$-a.e.\ $(y_1,\ldots,y_{m})$. 
Applying  \eqref{chaos} yields \eqref{678a}.

For the converse implication, we assume that \eqref{678a}
holds for $k\ge 1$. (The case $k=0$ follows from Lemma \ref{l98}.)
By the uniqueness part of Theorem \ref{tchaos0} we then
have $\|T_mf\|^2_m=0$ for all $m\ge k+1$ and hence
$$
\BE D^n_{y_1,\ldots,y_n}D^k_{x_1,\ldots,x_k}f(\eta)=0
\quad \lambda^n\text{-a.e.\ $(y_1,\ldots,y_n)$},\;
\lambda^k\text{-a.e.\ $(x_1,\ldots,x_k)$}
$$
for all $n\ge 1$. Applying \eqref{varcha} to the
function $D^k_{x_1,\ldots,x_k}f$ gives
\begin{align*}
\BE (D^k_{x_1,\ldots,x_k}f(\eta))^2=(\BE D^k_{x_1,\ldots,x_k}f(\eta))^2
\quad \lambda^k\text{-a.e.\ $(x_1,\ldots,x_k)$}.
\end{align*}
Since Jensen's inequality $\BE X^2\ge (\BE X)^2$ is an equality iff
$X$ is a.s.\ constant, we obtain that
\begin{align*}
D^k_{x_1,\ldots,x_k}f(\eta)=\BE D^k_{x_1,\ldots,x_k}f(\eta)
\quad \BP\text{-a.s.},\;\lambda^k\text{-a.e.\ $(x_1,\ldots,x_k)$}.
\end{align*}
By Lemma \ref{l98} and Fubini's theorem we get
$D^{k+1}_{x_1,\ldots,x_{k+1}}f(\eta)
=D^1_{x_{k+1}}D^{k}_{x_1,\ldots,x_{k}}f(\eta)=0$ 
almost surely for $\lambda^{k+1}$-a.e.\ $(x_1,\ldots,x_{k+1})$.

For the final part, assume \eqref{678a} is true. By 
the uniqueness part of
Theorem \ref{tchaos0} we  have ($\lambda^j$-a.e.)  
that $T_jf= f_j$ for $j \leq k$, as asserted. \qed

\vspace*{0.3cm}
Next we use Theorem \ref{tcov} to prove a series of
lower and upper variance bounds 
associated with the truncated series \eqref{varcha}.
The first of these upper bounds, i.e. the case $k=1$
of \eqref{Poin5} below, is the  Poincar\'e inequality
\eqref{poin1}
and, as remarked earlier, is also a special case of
Theorem \ref{tcovineq}. 

\begin{theorem}\label{tpoisecond} Let $f\in L^2(\BP_\eta)$
and $k\in\N$. Then 
\begin{align}\label{Poin5}
\sum^k_{n=1}\frac{1}{n!}
 \|T_nf\|^2_n \le 
\BV[f(\eta)]
\le \Big( \sum^{k-1}_{n=1}\frac{1}{n!}
 \|T_nf\|^2_n \Big) +\frac{1}{k!}\BE\| D^kf(\eta)\|^2_k,
\end{align}
with  an empty sum
 interpreted as zero.
Both inequalities are strict unless
$D^{k+1}_{y_1,\ldots,y_{k+1}}f(\eta) =0$ almost surely
for $\lambda^{k+1}$-a.e. $(y_1,\ldots,y_{k+1})$,  
in which case 
both inequalities are equalities.
\end{theorem} 
{\sc Proof:}  
The first inequality follows directly
from \eqref{varcha}. To prove the second, write
\eqref{varcha} as
\begin{align*}
\BV[f(\eta)] 
& = \Big( \sum^{k-1}_{n=1}\frac{1}{n!}
\|T_nf\|_n^2 \Big)
+\frac{1}{k!}\int (\BE D^k_{y_1,\ldots,y_k}
f(\eta))^2 \lambda^k(d(y_1,\ldots,y_k))\\
&+\sum^\infty_{n=k+1}\frac{1}{n!}
\iint (\BE D^{n-k}_{y_{k+1},\ldots,y_{n}}D^k_{y_1,\ldots,y_k}f(\eta))^2 
\lambda^{n-k}(d(y_{k+1},\ldots,y_{n}))\lambda^k(d(y_1,\ldots,y_k)).
\end{align*}
By relabelling $n-k$ as $m$ for $n>k$  and $y_\ell $ as $x_{\ell-k}$ for
$\ell > k$, 
it follows  that
\begin{align}
\notag
\BV[f(\eta)]
&\le \Big( \sum_{n=1}^{k-1} \frac{1}{n!} \|T_nf\|_n^2 \Big) + \frac{1}{k!}
 \int \Big[ (\BE D^k_{y_1,\ldots,y_k}f(\eta))^2 
\\
& +\sum_{m=1}^\infty\frac{1}{m!}
\int (\BE D^m_{x_1,\ldots,x_{m}}D^k_{y_1,\ldots,y_k}f(\eta))^2
 \lambda^m(d(x_1,\ldots,x_{m})) \Big]
\lambda^k(d(y_1,\ldots,y_k)).\notag
 \end{align}
Assuming without loss of generality that
$\BE\int (D^k_{y_1,\ldots,y_k}f(\eta))^2 \lambda^k(d(y_1,\ldots,y_k))<\infty$
(else the right hand side of \eqref{Poin5} is infinite),
we can apply
\eqref{varcha} to the function $D^k_{y_1, y_2,\ldots,y_k}f$
for $\lambda^k$-a.e.\ $(y_1,\ldots,y_k)$,
thereby simplifying the expression inside $[ \cdot ]$ above,
to obtain
\begin{align} \notag
\BV[f(\eta)]
&\leq \Big( \sum_{n=1}^{k-1}  \frac{1}{n!} \|T_nf\|_k^2 \Big) + \frac{1}{k!}
\int \BE (D^k_{y_1,\ldots,y_k}f(\eta))^2 
\lambda^k(d(y_1,\ldots,y_k)) 
\end{align}
which yields \eqref{Poin5}.

Assume now that one of the inequalities in \eqref{Poin5} is
an equality. Then $\|T_nf\|_n =0$ for all $n > k$. 
Therefore we obtain
from \eqref{chaos} and \eqref{orth} that \eqref{678a} holds
for  some $f_j\in \bH_j$, $j=1,\ldots,k$. 
Hence, by Proposition \ref{pfinchaos},
$D^{k+1}_{y_1,\ldots,y_{k+1}}f(\eta) =0$ almost surely
for $\lambda^{k+1}$-a.e.\ $(y_1,\ldots,y_{k+1})$,  

Assume, conversely, that
$D^{k+1}_{y_1,\ldots,y_{k+1}}f(\eta) =0$ almost surely
for $\lambda^{k+1}$-a.e. $(y_1,\ldots,y_{k+1})$.  
Then by Proposition \ref{pfinchaos}, \eqref{678a} holds 
for  some $f_j\in \bH_j$, $j=1,\ldots,k$. 
Therefore, by the uniqueness part of Theorem \ref{tchaos0} we  have
($\lambda^j$-a.e.)  that  $T_jf =0$ for $j > k$,
so that $\|T_nf\|_n =0$ for $n >k$ and hence      
we have equalities in  \eqref{Poin5}.\qed

\vspace*{0.3cm}
{\sc Proof of Theorem \ref{tcovineq}:} For any $n\in\N$ we abbreviate 
$a_n:=\|\BE D^nf(\eta)\|^2_n$.
Using \eqref{varcha} for $D^n_{y_1,\ldots,y_n}f$ we get
$$
\BE \int (D^n_{y_1,\ldots,y_n}f(\eta))^2\lambda^n(d(y_1,\ldots,y_n))
=\sum^\infty_{m=0}\frac{a_{m+n}}{m!}
=\sum^\infty_{m=n}\frac{a_{m}}{(m-n)!}.
$$
Hence, using \eqref{varcha} again, we have that  \eqref{covineq} 
is equivalent to
$$
\sum^{2k}_{n=1}\sum^\infty_{m=n}\frac{(-1)^{n+1}}{n!}\frac{a_{m}}{(m-n)!}
\le \sum^\infty_{m=1}\frac{a_{m}}{m!}
\le \sum^{2k-1}_{n=1}\sum^\infty_{m=n}\frac{(-1)^{n+1}}{n!}\frac{a_{m}}{(m-n)!}.
$$
Interchanging the order of summation shows that these inequalities
are implied by
$$
\sum^{(2k)\wedge m}_{n=1}\frac{(-1)^{n+1}}{n!}\frac{1}{(m-n)!}\le\frac{1}{m!}
\le \sum^{(2k-1)\wedge m}_{n=1}\frac{(-1)^{n+1}}{n!}\frac{1}{(m-n)!},\quad m\in\N,
$$
where $(2k)\wedge m$ denotes the minimum of $2k$ and $m$. Since the latter
equalities are equivalent to the elementary inequalities
\begin{align}\label{0811b}
\sum^{(2k)\wedge m}_{n=0}(-1)^{n}\binom{m}{n}\ge 0
\ge \sum^{(2k-1)\wedge m}_{n=0}(-1)^{n}\binom{m}{n},
\end{align}
this concludes the proof of \eqref{covineq}. Moreover, the first inequality
in \eqref{0811b} is strict unless $m \leq 2k$, so that the first
inequality of \eqref{covineq} is an equality if and only if
$a_m =0$ for all $m > 2k$. By Theorem \ref{tchaos0} 
and Proposition \ref{pfinchaos}, this
happens if and only if 
$D^{2k+1}_{y_1,\ldots,y_{2k+1}}f(\eta) =0$ almost surely
for $\lambda^{2k+1}$-a.e. $(y_1,\ldots,y_{2k+1})$.  

Similarly, the second  inequality
in \eqref{0811b} is strict unless $m \leq 2k -1 $, so that the second
inequality of \eqref{covineq} is an equality if and only if
$a_m =0$ for all $m \leq 2k-1$, which 
happens if and only if 
$D^{2k}_{y_1,\ldots,y_{2k}}f(\eta) =0$ almost surely
for $\lambda^{2k}$-a.e.\ $(y_1,\ldots,y_{2k})$.  
\qed

\section{Covariance identities}\label{scov}
\setcounter{equation}{0}

In this section we shall prove Theorem \ref{tkey} but first
we give some of its consequences.

\begin{theorem}\label{tkey2} 
Consider a Poisson process $\tilde\eta$ on $\tilde\BY:=\R_+\times\BY$
whose intensity measure is the product of Lebesgue measure
on $[0,1]$ and $\lambda$. 
Assume that $\eta=\tilde\eta([0,1]\times\cdot)$, and for $s > 0$
let $\tilde \eta_{s-}$
denote the restriction of $\tilde \eta $ to $[0,s) \times \BY$.
Let $f,g\in L^2(\BP_\eta)$. 
Then
\begin{align}\label{1b}
\BE\int\int^1_0
\BE[D_yf(\eta)|\tilde\eta_{s-}]^2ds\lambda(dy)<\infty
\end{align}
and
\begin{align}\label{Cov2}
\CV[f(\eta),g(\eta)]=\BE \int\int^1_0\BE[D_yf(\eta)|\tilde\eta_{s-}]
\BE[D_yg(\eta)|\tilde\eta_{s-}]ds\lambda(dy).
\end{align}
\end{theorem}
\noindent
{\sc Proof:} Apply Theorem \ref{tkey} to $\tilde \eta$ with the 
relation $<$ on $\tilde \BY$ given by $(s,x) < (t,y)$
if and only if $s < t$. 

\vspace{0.3cm}
\noindent
{\sc Proof of Theorem \ref{c2}:} 
It is no loss of generality
to assume that $\eta=\tilde\eta([0,1]\times\cdot)$, with $\tilde\eta$
as in Theorem \ref{tkey2}.
Then the  result is a direct
consequence of Theorem \ref{tkey2}.\qed

\begin{remark}\label{reiso}\rm
Again assuming without loss of generality
that $\eta=\tilde\eta([0,1]\times\cdot)$, 
taking $g =f$ and applying the conditional Jensen inequality
in \eqref{Cov2} yields an alternative proof of the Poincar\'e inequality
\eqref{poin1}. 
\end{remark}

We proceed with proving Theorem \ref{tkey}. 
For the rest of this section we assume that $\eta$ is a Poisson process on
$\BY = \R_+\times\BX$ with an intensity measure satisfying \eqref{diffuse}.
We need the following lemma.

\begin{lemma}\label{l97} Let $f\in L^2(\BP_\eta)$ and define
$f_{y}(\mu):=\int D_{y}f(\mu_{y}+\nu)\Pi^y(d\nu)$,
for $\mu\in \bN$, $y\in\BY$. Then 
\begin{align}\label{5.2}
\BE D^n_{(x_1,\ldots,x_n)}f_{y}(\eta)=
\Big( \prod_{i=1}^n  
\I\{x_i<y\} \Big)\BE D^{n+1}_{(x_1,\ldots,x_n,y)}f(\eta)
\end{align}
holds for $\lambda^{n+1}$-a.e.\ $(x_1,\ldots,x_n,y)$.
\end{lemma}
{\sc Proof:} 
For any $n \in \N$, and any $C \in \YY_0$,
arguing as in the proof of Lemma \ref{lemsubs} 
yields
$$
\int_{C^n} \BE \left[ \left|f\left(\eta + \sum_{i=1}^n
 \delta_{y_i}\right) \right| 
\right]
\lambda^n(d(y_1,\ldots,y_n)) < \infty.
$$ 
Hence
$\BE [ |f(\eta + \sum_{i=1}^n \delta_{y_i}) |] $ is finite
for  $\lambda^n$-almost all $(y_1,\ldots,y_n) \in \BY^n$. 
Thus for any $m \in \N_0$,
\begin{align}\label{0727a}
\int|f(\mu_{y}+\delta_{y_1}+\ldots+\delta_{y_m}
+\delta_{y}+\nu)|\Pi^y(d\nu)<\infty
\end{align}
and
\begin{align}\label{0730a}
\int|f(\mu_{y}+\delta_{y_1}+\ldots+\delta_{y_m }+\nu)|\Pi^y(d\nu)<\infty
\end{align}
hold for $\BP_\eta$-a.e.\ $\mu$ and
$\lambda^{m+1}$-a.e.\ $(y_1,\ldots,y_m,y)$,
since integrating the left hand side of
\eqref{0727a} over $\mu$ 
(with respect to the measure $\BP_\eta$) yields  
$\BE |f(\eta + \delta_{y} + \sum_{i=1}^{m}  \delta_{y_i} ) | $ and
integrating the left hand side of \eqref{0730a} over $\mu$ 
yields  $\BE |f(\eta +  \sum_{i=1}^{m}  \delta_{y_i} ) | $.

Let $x_i \in \BY$ for $1 \leq i \leq n$.
First suppose for some $i$ that $x_i < y$ does not hold.
Then since $f_{y}(\mu)$ depends on $\mu$ only through $\mu_{y}$,
and since $(\mu+\delta_{x_i})_{y} = \mu_{y}$, it follows
that $D_{x_i} f_{y}(\mu) =0$ for any $\mu$.
Since  $D_{y_1,\ldots,y_n}f(\mu)$ is symmetric
in $y_1,\ldots,y_n$ it follows that  
in this case \eqref{5.2} holds with both sides equal to zero.

Now suppose that $x_i < y$ for each $i$, and that \eqref{0727a}
and \eqref{0730a}
hold for every subset $\{y_1,\ldots,y_m\}$ of
$\{x_1,\ldots,x_n\}$. Then  by \eqref{Dsymmetric}
\begin{align}
D^n_{x_1,\ldots,x_n}f_{y}(\mu)&=
\sum_{J \subset \{1,\ldots,n\}} (-1)^{n-|J|} 
f_{y}\Big(\mu+\sum_{i\in J}\delta_{x_i}\Big)\notag \\
&= \sum_{J \subset \{1,\ldots,n\}} (-1)^{n-|J|} 
\int D_{y} f \Big(\Big(
\mu + \sum_{i\in J}\delta_{x_i} \Big)_{y} + \nu \Big) \Pi^y(d \nu)
\notag  \\
&=\int D^{n+1}_{y,x_1,\ldots,x_n} f (\mu_{y} + \nu) \Pi^y(d \nu).
\label{419}
\end{align}
Integrating over $\mu$ (with respect to the measure $\BP_\eta$)
yields \eqref{5.2} for this case. \qed

\vspace{0.3cm}
{\sc Proof of Theorem \ref{tkey}:} Let us first assume that
$f,g$ are bounded. Then $f_{y}$ and $g_{y}$ are trivially
in $L^2(\BP_\eta)$ for all $y\in\BY$ and we have
by \eqref{956} that
\begin{align*}
\BE\int \BE[D_{y}f(\eta)|\eta_{y}]
\BE[D_{y}g(\eta)|\eta_{y}]\lambda(dy)
=\int \BE[f_{y}(\eta) g_{y}(\eta)]\lambda(dy),
\end{align*}
where the use of Fubini's theorem will be justified below.
By Theorem \ref{tcov} and Lemma \ref{l97} this equals
\begin{align*}
\int\sum_{n=0}^\infty\frac{1}{n!}
\langle T_nf_{y},T_ng_{y}\rangle_{n}\lambda(dy)
=\sum_{n=0}^\infty\frac{1}{n!}
\langle\I_{J_{n+1}}T_{n+1}f,T_{n+1}g\rangle_{n+1},
\end{align*}
where $J_{n+1}:=\{(y_1,\ldots,y_{n+1}):y_1<
y_{n+1},\ldots,y_n<y_{n+1}\}$
and where the interchange of summation 
and integration will be justified below.
By \eqref{diffuse} the measure $\lambda^{n+1}$ is concentrated on 
$\cup_{\pi \in 
\Sigma_{n+1}}\{(y_1,\ldots,y_{n+1}):y_{\pi(1)}<\ldots<y_{\pi(n+1)}\}$
where the union is over all permutations.
Using the symmetry of $T_{n+1}f$ and $T_{n+1}g$ this gives
$$
\frac{1}{n!}\langle\I_{J_{n+1}}T_{n+1}f,T_{n+1}g\rangle_{n+1}
=\frac{1}{(n+1)!}\langle T_{n+1}f,T_{n+1}g\rangle_{n+1}.
$$
Using Theorem \ref{tcov} again, yields the asserted identity \eqref{Cov}.
Repeating the above calculation with $f=g$ shows that \eqref{sqint}
holds. This  also justifies our use of Fubini's theorem.

Finally we consider the case of general $f,g\in L^2(\BP_\eta)$. 
The previous arguments carry through once we have shown that
\eqref{sqint} holds for $f\in L^2(\BP_\eta)$. Let $f^k$, $k\in\N$, be a sequence
of bounded measurable functions on $\bN(\BY)$ such
$\BE (f(\eta)-f^k(\eta))^2\to 0$ as $k\to\infty$.  We have just proved
that
\begin{align*}
\BV [f^k(\eta)-f^l(\eta)]=
\BE\int(\BE[D_{y}f^k(\eta)|\eta_{y}]-\BE[D_{y}f^l(\eta)
|\eta_{y}])^2\lambda(dy), 
\quad k,l\in\bN.
\end{align*}
Since $L^2(\BP\otimes\lambda)$ is complete,
there is an $h\in L^2(\BP\otimes\lambda)$ satisfying  
\begin{align}\label{345}
\lim_{k\to\infty}\BE\int(h(y)-\BE[D_{y}f^k(\eta)|\eta_{y}])^2\lambda(dy)=0.
\end{align}
On the other hand it follows from Lemma \ref{lemsubs} that
for any $C\in\mathcal{Y}_0$
\begin{align*}
\int_C \BE\big|\BE[D_{y}f^k(\eta)|\eta_{y}]&-\BE[D_{y}f(\eta)
|\eta_{y}]\big|\lambda(dy)\\ 
&\le \int_C \BE |D_{y}f^k(\eta)-D_{y}f(\eta)|\lambda(dy)\to 0
\end{align*}
as $k\to\infty$. Comparing this with \eqref{345} shows that
$h(\omega,y)=\BE[D_{y}f|\eta_{y}](\omega)$ for 
$\BP\otimes\lambda$-a.e.\ $(\omega,y)\in\Omega\times C$
and hence also for $\BP\otimes\lambda$-a.e.\ $(\omega,y)\in\Omega\times \BY$.
Therefore the fact that
$h\in L^2(\BP_\eta\otimes\lambda)$ implies \eqref{sqint}.\qed

\section{Infinitely divisible random measures}\label{sinfdiv}
\setcounter{equation}{0}

In this section we consider an infinitely divisible random measure $\xi$
on a complete separable metric space $\BX$
equipped with the Borel $\sigma$-field $\mathcal{X}$, see
\cite{MKM,DVJ} and \cite{Ka83} for the the special case
of a locally compact phase space $\BX$. The distribution
of such a random measure can be most conveniently described
by its {\em Laplace functional} as follows.
Let $\bM$ denote the space of all locally finite measure on $\BX$, that is,
the set of all measures that are finite
on metrically bounded sets. We equip $\bM$ with
the smallest $\sigma$-field of subsets of $\bM$ 
such that the mappings $\mu\mapsto\mu(B)$ are measurable for all
$B\in\mathcal{X}$.
Let $v:\BX\rightarrow\R_+$ be measurable. Then 
\begin{align}\label{Lap}
\BE \exp\left[-\int v d\xi \right]=
\exp\left[-\int v d\alpha-\int \big(1-e^{-\int vd\mu}\big)\BQ(d\mu)\right],
\end{align}
where $\alpha\in\bM$ and the {\em KLM} 
(or {\em L\'evy}) {\em measure} $\BQ$ is a $\sigma$-finite measure
on $\bM$ having $\BQ(\{0\})=0$ (here $0$ denotes the
zero measure) and 
\begin{align}\label{fin}
\int(1-\exp[-\mu(B)])\BQ(d\mu)<\infty, \quad B\in\mathcal{X}_0.
\end{align}
Here $\mathcal{X}_0\subset\mathcal{X}$ denotes the ring
of (metrically) bounded Borel subsets of $\BX$.

\begin{proposition}\label{tP}
For any $f\in L^2(\BP_{\xi})$,
\begin{align}\label{Poin2}
\BV[f(\xi)]\le \BE\int (f(\xi+\mu)-f(\xi))^2\BQ(d\mu).
\end{align}
\end{proposition}
{\sc Proof:} Let $C_n$, $n\in\N$, 
be a sequence of closed balls
with fixed centre and radius $n$. We define a measurable
mapping $H:\bN(\bM)\rightarrow\bM$ as follows.
Let $\nu\in \bN(\bM)$ and define
\begin{align*}
H(\nu)(B):=\alpha(B)+\int \mu(B)\nu(d\mu),\quad B\in\mathcal{X},
\end{align*}
whenever the right hand side is finite for all $B=C_n$,
$n\in\N$. Otherwise we define $H(\nu):=0$. 
Now let $\eta$ be a Poisson process on $\bM$ with
intensity measure $\BQ$. It is well-known \cite{Ka83,MKM,DVJ}
that $\xi$ and $H(\eta)$ have the same distribution.
Because of this {\em Poisson cluster representation} of $\xi$
we can assume $\xi=H(\eta)$ without restricting
generality. Now we apply Theorem \ref{tcovineq} to the
function $f\circ H$. Since $H(\nu+\delta_\mu)=H(\nu)+\mu$
for all $\nu\in\bN(\bM)$ and all $\mu\in\bM$
we obtain the assertion.\qed

\vspace*{0.3cm}
In case $\alpha=0$ and $\BQ=\int\I\{\delta_x\in\cdot\}\lambda'(dx)$
for some $\lambda'\in \bM$ 
the random measure $\xi$ is a Poisson process with
intensity measure $\lambda'$. Then
the inequality \eqref{Poin2} simplifies to the Poincar\'e inequality 
\eqref{poin1}. 

By means of the Poisson cluster representation of $\xi$
we could also rewrite the other results of this paper.
We restrict ourselves to the following version of 
the Harris-FKG inequality.
We call a measurable function $f:\bM\rightarrow\R$ 
{\em increasing almost everywhere} if
$f(\xi+\mu)\ge f(\xi)$ $\BP$-a.s.\ and for $\BQ$-a.e.\ $\mu$.
The next proposition follows from Theorem \ref{c2}. 

\begin{proposition}\label{tFKG}
Assume that $f,g\in L^2(\BP_{\xi})$ are increasing almost everywhere.
Then
\begin{align}\label{FKG2}
\BE [f(\xi)g(\xi)]\ge (\BE f(\xi))(\BE g(\xi)).
\end{align}
\end{proposition}

\section{The case of finite $\BY$}\label{sfinite}
\setcounter{equation}{0}

For this section only, we assume $\BY$ is a finite set $\{1,\ldots,k\}$
and $\YY$ is the power set of $\BY$. In this case let us write
$\lambda_i$ for $\lambda(\{i\})$ (assumed finite)
and $\eta_i $ for $\eta(\{i\})$.
Then $\eta_1, \ldots,\eta_k$ are independent Poisson variables.
Given $n \in \N$ and $\by = (y_1,\ldots,y_n) \in \BY^n$,
let $\bm(\by) = (m_1(\by),\ldots,m_k(\by))$ be given  
by
$$
m_j(\by) := \sum_{i=1}^n {\bf 1} \{ y_i =j\},
$$
so that $\sum_{j=1}^k m_j(\by) = n$.
Then for any $h: \BY^n \to \R$,
by \eqref{stcharl} and linearity, 
\begin{align}\label{0730d}
I_n(h)=\sum_{\by \in \BY^n} h(\by) \prod_{i=1}^k 
\Big(\lambda_i^{m_i} C_{m_i}(\lambda_i;\eta_i)\Big). 
\end{align}
Now consider $f:\N_0^k\to\R$ with $\BE f(\eta_1,\ldots,\eta_k)^2< \infty$.
Then  $T_nf(y_1,\ldots,y_n)$ depends on $\by$ only through $\bm := \bm(\by)$. 
Writing $m_i$ for $m_i(\by)$,  and using \eqref{Dsymmetric},  we define
$$
a_{\bm} : = T_n f(\by) = \sum_{i_1=0}^{m_1} \ldots \sum_{i_k=0}^{m_k} 
(-1)^{n - \sum_{j=1}^k i_j} \left( \prod_{j=1}^k \binom{m_j}{i_j} \right) 
\BE [ f(\eta_1 +i_1,\ldots,\eta_k +i_k) ].
$$
Then by \eqref{0730d},
$$
I_n(T_nf) = \sum_{\bm: m_1 + \ldots + m_k = n} \Big(
\frac{n!}{m_1! \cdots m_k!} \Big) 
a_{\bm} \prod_{i=1}^k \Big( \lambda_i^{m_i} C_{m_i}(\lambda_i;\eta_i) \Big)
$$
and therefore by Theorem \ref{tchaos0},
\begin{align}\label{0730b}
f(\eta_1,\ldots,\eta_k) = 
\BE [ f(\eta_1,\ldots,\eta_k) ] + \sum_{\bm} a_{\bm} 
\prod_{i=1}^k \Big(\lambda_i^{m_i}C_{m_i}(\lambda_i;\eta_i)/m_i! \Big)
\end{align}
where the sum is over all $\bm = (m_1,\ldots,m_k) \in \N_0^k$ 
except for $(0,\ldots,0)$, and the convergence is in $L^2(\BP)$.

This identifies the coefficients in the Charlier polynomial
expansion of $f$ in terms of the expected repeated differences
$a_\bm$. Note that 
$\BE [ C_m(\lambda_i;\eta_i)^2] = m! \lambda_i^{-m}$
(to see this, use \eqref{0730d} and \eqref{orth} in
the case with $k=1$ and $h=g =1$). Hence by taking inner products with
$\prod_{i=1}^k C_{m_i}(\lambda_i;\eta_i)$ in \eqref{0730b} and using
orthogonality of the Charlier polynomials, we obtain
$$
a_{\bm} =
\BE \Big[  f(\eta_1,\ldots,\eta_k) \prod_{i=1}^k C_{m_i}(\lambda_i;\eta_i) \Big].
$$
In the case $k=1$ this has previously been obtained as 
Lemma 9.1.4 of \cite{BHJ92}. 


\vspace{0.3cm}
\noindent
{\bf Acknowledgements:} We wish to thank Yuri Kabanov for giving
several insightful comments on an early draft of
this paper and Liming Wu for sending us his unpublished
work \cite{Wu98}. Thanks are also due to Matthias Reitzner,
who drew our attention to the Poincar\'e inequality 
for general Poisson processes.


\begin{thebibliography}{99}

\bibitem{AneLed00}
An\'e, C. and Ledoux, M. (2000). 
On logarithmic Sobolev inequalities for continuous time random 
walks on graphs. {\it Probab.\ Theory Related Fields} {\bf 116}, 573-602.

\bibitem{BHJ92}
Barbour, A. D., Holst, L. and Janson, S. (1992). {\it Poisson Approximation.}
Clarendon Press, Oxford.

\bibitem{Bl95}
Blaszczyszyn, B. (1995). Factorial-moment expansion for stochastic systems. 
{\em Stoch. Proc. Appl.} {\bf 56}, 321-335. 

\bibitem{BR06}
Bollob\'as, B. and  Riordan, O.
(2006).
{\em Percolation.} Cambridge University Press, New York, 

\bibitem{Chen85}
Chen, L.\ (1985). Poincar\'e-type inequalities via stochastic integrals. 
{\em Z.\ Wahrscheinlichkeitstheorie verw. Gebiete} {\bf 69}, 251-277.


\bibitem{DVJ}
Daley, D.J.\ and Vere-Jones, D.\ (2008).
\newblock {\it An Introduction to the Theory of Point Processes, Volume II}.
2nd ed., 
\newblock Springer, New York.

\bibitem{DeMeyer78}
Dellacherie, C.\ and Meyer, P.A.\ (1978).
Probabilities and potential. Mathematics Studies, Volume 29, 
Hermann, Paris. 
North-Holland Publishing Company, Amsterdam and New York.




\bibitem{ES}
Efron, B. and Stein, C. (1981).
The jackknife estimate of variance.
{\em Ann. Statist.} {\bf 9},  586--596.



\bibitem{GeK97}
Georgii, H.\ and K\"uneth, T.\ (1997).
Stochastic order of point processes. 
{\em J. Appl. Prob.} {\bf 34}, 868-881. 

\bibitem{HevRei}
Heveling, M. and   Reitzner, M. (2009).
Poisson-Voronoi approximation
{\em Ann. Appl. Probab.} {\bf 19},  719-736.

\bibitem{Hitsuda72}
Hitsuda, M.\ (1972). Formula for Brownian partial derivatives. 
Proceedings of the 2nd Japan-USSR Symposium on Probability Theory. 
111-114. 


\bibitem{HoPer95}
Houdr\'e, C., and Perez-Abreu, V. (1995). 
Covariance identities and inequalities for 
functionals on Wiener space and Poisson space. 
{\em Ann.\ Probab.} {\bf 23}, 400-419.

\bibitem{Ho02}
Houdr\'e, C.\ (2002).
Remarks on deviation inequalities for functions of 
infinitely divisible random vectors.
{\it Ann.\ Probab.} {\bf 30}, 1223-1237. 

\bibitem{HoP02}
Houdr\'e, C.\ and Privault, N.\ (2002).
Concentration and deviation inequalities in 
infinite dimensions via covariance representations.
{\em Bernoulli} {\bf 8}, 697-720. 

\bibitem{Ito51}
It\^o, K. (1951).
Multiple Wiener integral. {\it J. Math. Soc. Japan} {\bf 3}, 157-169.

\bibitem{Ito56}
It\^o, K. (1956).
Spectral type of the shift transformation 
of differential processes with stationary 
increments. {\it Trans.\ Amer.\ Math.\ Soc.} {\bf 81}, 253-263.

\bibitem{Ito88}
Ito, Y.\ (1988). Generalized Poisson functionals. 
{\em Probab.\ Th.\ Rel.\ Fields} {\bf 77}, 1-28.


\bibitem{Kab75}
Kabanov, Y.M.\ (1975). On extended stochastic integrals. 
{\em Theory Probab. Appl.} {\bf 20}, 710-722.

\bibitem{KabSk75}
Kabanov, Y.M.\  and Skorokhod, A.V.\ (1975). 
Extended stochastic integrals. Proceedings of the School-Seminar on
the Theory of Random Processes. Druskininkai, November 25-30, 1974.
Part I. Vilnius, (Russian).


\bibitem{Ka83}
Kallenberg, O.\ (1983). 
\newblock {\it Random Measures.}  
\newblock Akademie-Verlag Berlin and Academic Press, London.


\bibitem{Kallenberg}
Kallenberg, O.\ (2002). 
\newblock {\it Foundations of Modern Probability}. 
\newblock Second Edition, Springer, New York.



\bibitem{Lieb94}
Liebscher, V.\ (1994):
On the isomorphism of Poisson space and symmetric Fock Space.  
In L. Accardi, editor, {\em Quantum Probability \& Related Topics IX},
pp.\ 295-300, World Scientific, Singapore, .

\bibitem{Lo05}
L{\o}kka, A.\ (2005).  
Martingale representation of functionals of L\'evy processes.
{\it Stochastic Analysis and Applications} {\bf 22}, 867-892.

\bibitem{MKM}
Matthes, K., Kerstan, J.\ and Mecke, J.\ (1978).
{\it Infinitely Divisible Point Processes}.
Wiley, Chichester.


\bibitem{Mecke} 
Mecke, J.\ (1967). 
\newblock{Station\"are zuf\"allige Ma\ss e auf lokalkompakten Abelschen Gruppen}. 
{\em Z. Wahrsch. verw. Gebiete} {\bf 9}, 36-58.


\bibitem{MR96}
Meester, R. and Roy, R.\ (1996).
{\it Continuum Percolation}, Cambridge Univ.\ Press. 

\bibitem{Meyer95}
Meyer, P.A.\ (1995). {\em Quantum Probability for Probabilists}
2nd ed., Springer, Berlin. 

\bibitem{MZ960}
M{\o}ller, J.\ and Zuyev, S.\ (1996).
Gamma-type results and other related properties of Poisson processes.
{\em Adv.\ Appl.\ Prob.} {\bf 28}, 662-673.

\bibitem{MolZu00}
Molchanov, I.\ and Zuyev, S (2000).
Variational analysis of functionals of Poisson processes. 
{\em Math. Operat. Res.} {\bf 25}, 485-508.

\bibitem{NuViv90}
Nualart, D.\ and Vives, J.\ (1990). Anticipative calculus for the 
Poisson process based on the Fock space. 
In: S\'eminaire de Probabilit\'es XXIV, Lecture Notes in Math., 
{\bf 1426}, 154-165.


\bibitem{Ogura72}
Ogura, H.\ (1972).
Orthogonal functionals of the Poisson processes. 
{\it Trans IEEE Inf. Theory} {\bf 18}, 473-481.


\bibitem{PSTU09}
Peccati, G., Sol\'e, J. L., Taqqu, M.S.\  and Utzet, F.\ (2009).
Stein's method and normal approximation of Poisson functionals.
{\em Ann. Probab.}, to appear.

\bibitem{Pen03}
Penrose, M.\ (2003).
{\it Random geometric graphs}. Oxford University Press, Oxford.

\bibitem{PS05}
Penrose, M. D. and Sudbury, A.  (2005).
Exact and approximate results for deposition and annihilation processes 
on graphs. 
{\em Ann. Appl. Probab.} {\bf 15}
 853--889. 

\bibitem{PW08}
Penrose, M.D.\ and Wade, A.R.\ (2008). Multivariate normal
approximation in geometric probability.
{\em Journal of Statistical Theory and Practice} {\bf 2},
293-326.


\bibitem{Pic96}
Picard, J.\ (1996). Formules de dualit\'e sur l’espace de Poisson.
{\it Ann.\ Inst.\ H.\ Poincar\'e  Probab.\ Statist.} {\bf 32}, 509-548. 

\bibitem{Pic96b}
Picard, J.\ (1996). On the existence of smooth densities for jump processes. 
{\it Probab.\ Theory Related Fields} {\bf 105}, 481-511. 

\bibitem{Priv01}
Privault, N. (2001).
Extended covariance identities and inequalities. 
Statistics \& Probability Letters {\bf 55},  247-255.


\bibitem{Surg84}
Surgailis, D. (1984)
On multiple Poisson stochastic integrals and associated Markov semigroups.
{\it Probab. Math. Statist.} {\bf 3}, 217-239. 

\bibitem{Skor75}
Skorohod, A. V.\ (1975). On a generalization of a stochastic integral. 
{\it Theory Probab. Appl.} {\bf 20}, 219-233. 

\bibitem{Wiener38}
Wiener, N.\ (1938).
The homogeneous chaos. {\it Am J. Math.} {\bf 60}, 897-936.

\bibitem{Wu98}
Wu, L.\ (1998).  $L^1$ and modified logarithmic Sobolev inequalities 
and deviation inequalities for Poisson point processes. Preprint.

\bibitem{Wu00}
Wu, L.\ (2000).
A new modified logarithmic Sobolev inequality for 
Poisson point processes and several applications.
{\it Probab.\ Theory Related Fields} {\bf 118}, 427-438.


\end{thebibliography}
\end{document}